\documentclass[journal]{IEEEtran}
\usepackage[cmex10]{amsmath}
\usepackage[cmex10]{amsmath}

\usepackage{ifpdf}

\usepackage{algorithm}
\usepackage{algpseudocode}
\usepackage{color}
\usepackage{url}
\ifpdf
\usepackage[pdftex]{graphicx}
\else
\usepackage[dvips]{graphicx}
\fi
\usepackage{tabularx}
\usepackage{multirow}
\graphicspath{{graphics/png/}}
\DeclareGraphicsExtensions{.png, .eps, .pdf}
\usepackage{float} 
\usepackage{subfigure}
\usepackage[caption=false,font=footnotesize]{subfig}
\usepackage{comment}
\usepackage{bm}
\usepackage{amsmath,amsthm,amsfonts,amssymb}
\usepackage[noadjust]{cite}
\usepackage{nomencl}
\theoremstyle{definition}
\newtheorem{theorem}{\textbf{Theorem}}

\newtheorem{remark}[theorem]{\textbf{Remark}}

\newcommand{\up}[1]{^\mathrm{#1}}
\newcommand{\rev}[1]{\textcolor{black}{#1}}

\IEEEoverridecommandlockouts
\begin{document}
%
% paper title
% Titles are generally capitalized except for words such as a, an, and, as,
% at, but, by, for, in, nor, of, on, or, the, to and up, which are usually
% not capitalized unless they are the first or last word of the title.
% Linebreaks \\ can be used within to get better formatting as desired.
% Do not put math or special symbols in the title.
\title{Arbitraging Variable Efficiency Energy Storage using Analytical Stochastic Dynamic Programming}

%% To specify the authors when (number of affiliations <= 2)
\author{Ningkun Zheng,~\IEEEmembership{Student Member,~IEEE,} Joshua Jaworski, Bolun~Xu,~\IEEEmembership{Member,~IEEE}
	\thanks{N.~Zheng, J. Jaworski, and B.~Xu are with Columbia University, NY, USA.}
}

%% To specify the authors when (number of affiliations > 2)
% \author{\IEEEauthorblockN{Author n.1\IEEEauthorrefmark{1},
% Author n.2\IEEEauthorrefmark{2},
% Author n.3\IEEEauthorrefmark{3}, 
% Author n.4\IEEEauthorrefmark{3} and
% Author n.5\IEEEauthorrefmark{4}}
% \IEEEauthorblockA{\IEEEauthorrefmark{1} Department Name of Organization A\\
% Name of the organization A,
% Address A\\ Emails if wanted}
% \IEEEauthorblockA{\IEEEauthorrefmark{2} Department Name of Organization B\\
% Name of the organization B,
% Address B\\ Emails if wanted}
% \IEEEauthorblockA{\IEEEauthorrefmark{3} Department Name of Organization C\\
% Name of the organization C,
% Address C\\ Emails if wanted}
% \IEEEauthorblockA{\IEEEauthorrefmark{4}Department Name of Organization D\\
% Name of the organization D,
% Address D\\ Emails if wanted}
% }

% make the title area
\maketitle

% As a general rule, do not put math, special symbols or citations
% in the abstract
\begin{abstract}
This paper presents a computation-efficient stochastic dynamic programming algorithm for solving energy storage price arbitrage considering variable charge and discharge efficiencies. We formulate the price arbitrage problem using stochastic dynamic programming and model real-time prices as a Markov process. Then we propose an analytical solution algorithm using a piecewise linear approximation of the value-to-go function.  Our solution algorithm achieves extreme computation performance and solves the proposed arbitrage problem for one operating day in less than one second on a personal computer. We demonstrate our approach using historical price data from four  price zones in New York Independent System Operator, with case studies comparing the performance of different stochastic models and storage settings. Our results show that the proposed method captures 50\% to 90\% of arbitrage profit compared to perfect price forecasts.  In particular, our method captures more than 80\% of arbitrage profit in three out of the four price zones when considering batteries with more than two-hour duration and realistic degradation cost.

% Instead of sampling discrete scenarios, our method calculates the storage value function from real-time price predicted by Markov process model. The proposed approach provides a low computational cost analytical approach to solve storage arbitrage , and case studies are performed on 4 different zones in NYISO.
\end{abstract}


\IEEEpeerreviewmaketitle
% Use this to place sponsorships
% \thanksto{Applicable sponsors, if any, should be placed using the \emph{thanksto} command.}

\section{Introduction}

Decarbonization in the electric power sector is pivotal for sustainable developments as other sectors become electrified. The renewable energy share of the US electricity generation mix has reached more than 20\% by 2020 and is projected will double by 2050~\cite{us2021annual}.  Electrical energy storage resources provide reliable and affordable solutions to many challenges in a power system with a high renewable share, such as peak shaving~\cite{oudalov2007sizing}, voltage support~\cite{krata2018real}, and frequency regulation~\cite{stroe2016operation}. Pilot programs have shown storage is a cost-effective and profitable solution in these applications. In a competitive market, storage that maximizes profit will also help to maximize social welfare~\cite{castillo2013profit}. The Federal Energy Regulatory Commission (FERC) issued Order 841~\cite{federal2018electric} to further remove barriers for energy storage to participate in wholesale electricity markets and urge independent system operators and regional transmission organizations to develop bidding and dispatch models for energy storage.

Yet, the participation of energy storage in wholesale energy markets has been limited compared to other applications, even with the fast dropping cost of energy storage~\cite{schmidt2017future}. As the most significant market and foundation of deregulated power systems, it is critical for storage to participate in wholesale energy markets efficiently. However, although all system operators have allowed energy storage to participate, current market designs have difficulties in dispatching storage accurately, restricted by the capability of forecasting and the limitation of computation power~\cite{martinez2020forecast}. Therefore, practical and efficient bid and control designs are heavily reliant on energy storage to better participate in energy markets.

This paper proposes an analytical stochastic dynamic programming (SDP) algorithm for optimizing variable efficiency energy storage price arbitrage in real-time energy markets with extreme computation efficiency. Our method targets a generic energy storage model with variable efficiency and discharge cost. Compared to optimization-based storage bidding and control methods such as bi-level optimization~\cite{mohsenian2015coordinated,pandvzic2015energy,wang2017look}, our method is lightweight and easy to implement. While compared to other dynamic programming or reinforcement learning methods~\cite{wang2018energy, jiang2015optimal}, our method does not discretize energy storage controls and states and models variable efficiency without deteriorating computation efficiency. Our method is efficient to train and can explicitly incorporate physical storage parameters such as efficiency and degradation. 
% Different approaches on energy storage arbitrage have been proposed in previous works. One category of approach is to use bi-level programming to optimize energy storage in day-ahead market~\cite{mohsenian2015coordinated,pandvzic2015energy,wang2017look}, but bi-level approaches are computation expensive and inefficient to incorporate uncertainty models or more sophisticated bidding settings. ~\cite{shafiee2016risk} proposed a risk-constrained operation scheduling for energy storage based on information gap decision theory to optimize energy storage bidding. ~\cite{wang2018energy} used a data-driven reinforcement learning algorithm for optimizing energy storage arbitrage. The main weakness of these study are over-simplify and discrete the charge/discharge actions, which made it unable to control ESRs accurately, or do not take account of stage-wise dependency of states. Our approach utilize a low computational power demanding stochastic dynamic programming (SDP) algorithm, while consider uncertainty in real-time price by integrating Markov process price prediction model.
The main contribution of this paper is as follows.

\begin{enumerate}
    \item We develop an analytical algorithm for solving storage SDP for real-time energy price arbitrage without using third-party solvers.  The algorithm uses less than one second to solve a real-time arbitrage problem for one operating day, including 288 stages (5-minute market clearing).
    \item We model real-time price uncertainties as a Markov process. The Markov models are trained using historical price data and are then incorporated into SDP to model price node transition probabilities.
    \item We compare different Markov process models for real-time price arbitrage, including stage and season dependency. The result shows that considering stage-wise dependency in price prediction improves model performance significantly.
    \item We demonstrate our algorithm by performing price arbitrage in different price zones in New York Independent System Operator (NYISO). Results show that the proposed method captures 50\% to 90\% of arbitrage profit compared to perfect price forecasts depending on the storage setting and price zones.
    \rev{\item We implement variable efficiency model to our proposed stochastic dynamic programming algorithm and benchmarked it with a piece-wise linear mixed-integer linear programming model (Julia+Gurobi) over deterministic arbitrage test cases.}
\end{enumerate}

The remainder of this paper is organized as follows. Section II presents the background and discusses previous works on this topic. Sections III formulates the problem, and Sections IV provides the solution method to proposed models. Section V includes simulation results and discusses the result compared to different benchmarks. Section VI concludes the paper and proposes a future course of research.

\section{Backgrounds and Literature Review}

% \subsection{Background}

% The increasing renewable penetration has prompted the deployment of energy storage to confront concerns regarding the reliability of power system~\cite{paul2010role}. Grid-scale energy storage in the US is projected to grow five-fold by 2050~\cite{nrel2021grid}. The Federal Energy Regulatory Commission (FERC) issued Order 841~\cite{ferc841} to further remove barriers for energy storage to participate in wholesale electricity market, and urge independent system operators/regional transmission organizations (ISOs/RTOs) to develop bidding and dispatch model for energy storage. Besides, FERC's recent Order 2222 allows distributed energy resources (DERs) and electric vehicles (EVs) to aggregate to participate in wholesale market, which also an edge implementation of energy storage arbitrage algorithm.

% \subsection{Literature Review}

Increasing penetration of renewable energy in power systems increases the fluctuation of electricity prices, and real-time market price arbitrage will become more profitable~\cite{woo2011impact}. Electrochemical battery energy storage can switch between full charging and discharging power in less than a second, providing it with a unique advantage to arbitrage real-time price differences~\cite{gevorgian2013ramping}. FERC Order 841 further ensures energy storage resources have equal access to wholesale electricity markets, and a growing number of researchers are focusing on battery energy storage arbitrage in the real-time markets.

Model predictive control (MPC) is one of the most widely used control methods in energy storage applications~\cite{arnold2011model}. MPC is an optimization-based control strategy. It solves an optimization problem over a finite future horizon with predictions, but only applies the first control step. MPC repeatedly optimizes the best-predicted system performance over a finite horizon with new system states update at each time step~\cite{maciejowski2002predictive}. Although MPC relies on some heuristic and expert knowledge in specific tasks, it achieves good performance in practice~\cite{morari1999model}. 

Due to MPC's excellent flexibility to incorporate different objectives and constraints,  previous works utilized MPC based controllers in energy storage applications. Shan et al.~\cite{shan2018universal} investigated energy storage for smoothing wind and solar power generation with MPC controllers. Khalid et al.~\cite{khalid2010model} presented a controller based on MPC for the efficient operation of the BESS for frequency regulation. Meng et al.~\cite{meng2015cooperation} proposed a distributed MPC scheme, which coordinates BESS for voltage regulation. The University of Queensland used an MPC controller for performing price arbitrage in electricity spot markets and achieved 70\% profit compared to the perfect forecast~\cite{qs}.  However, MPC must use a separate module to generate future scenarios and can be inefficient to incorporate uncertainty considerations. MPC is also computationally expensive, especially in edge implementations, because it has to repetitively solve optimization problems.

Stochastic dynamic programming (SDP) is a systematic approach to optimize price arbitrage in highly stochastic scenarios compared to MPC. SDP has been a classic approach to optimize hydropower and pump hydro scheduling but uses relatively simple stochastic models focused on weekly to seasonal scheduling~\cite{abgottspon2009mid, gjelsvik2010long,korpaas2015norwegian}. The major challenge of applying SDP in more sophisticated scenarios like price arbitrage is the computational cost, which grows exponentially as the number of decision stages and scenarios increases. Therefore, SDP has been less attractive in BESS real-time arbitrage due to the high frequent decision making nature (usually 5-minute) and prolonged decision-making horizon.

Nevertheless, a few exploratory attempts have been made to apply SDP in price arbitrage. One approach is to solve SDP using a Markov decision process (MDP), in which the storage control pairs with system and uncertainty states. It discretizes energy storage state of charge (SoC) and power, and has been applied in applications such as price arbitrage and frequency regulation~\cite{donadee2013optimal,wen2019optimal}. MDP has also been combined with reinforcement learning in energy storage price arbitrage~\cite{yu2020energy}, or used to investigate the welfare optimization considering consumer and producer surplus~\cite{deng2020community}.
However, these studies have to discretize states to map control to observed states. It induces a trade-off problem: the problem is too big to solve in reasonable times if states are discreitized to fine granularity, while coarse state granularity may not be accurate enough for storage arbitrage because of charging and discharging efficiency.

Stochastic dual dynamic programming (SDDP) is an alternative approach to solve SDP by approximating the value-to-go function using linear cuts based on the dual variable~\cite{shapiro2011analysis}. SDDP does not discretize storage control or states, and the computation performance is better than solving SDP using MDP~\cite{megel2015stochastic}. However, SDDP uses forward-backward iterations to generate dual cuts, and the number of iterations increases significantly with the number of uncertainty stages and nodes~\cite{aaslid2019optimal}. Comparably, the method proposed in this paper follows a similar theory foundation with SDDP, but we use an analytical algorithm to generate complete piecewise linear value function approximations. Our method solves the arbitrage SDP problem in a backward path without iterating. providing improved accuracy and computation speed.

\rev{Our proposed method also incorporates variable energy storage efficiency in arbitrage optimization without deteriorating computation performance.
Most energy storage technologies have variable efficiencies, lithium-ion batteries have near linear efficiencies, but other technologies such as flow batteries and hydrogen storage have strong nonlinearities~\cite{ibrahim2008energy}. Prior studies use mixed-integer piece-wise linearization to model nonlinear and non-convex efficiencies over a deterministic optimization, but few have explored incorporating variable efficiency with stochastic storage control primarily due to computation complexities~\cite{jafari2019improved,flamm2019price}. 
}
% In this paper, we use a piece-wise linearization approach of value function instead of doing state or control discretization, and employs an analytical algorithm to update value function dramatically improves computation speed. Furthermore, we consider step-wise dependency of real-time price by utilizing Markov process model to predict real-time price.

\section{Formulation}

% We start with formulating the energy storage arbitrage  as a stochastic dynamic programming problem, followed by introducing the Markov process used to model the stochastic price process. 
% In this section, we further extend previous work using State of Charge (SoC) valuation in energy storage arbitrage model~\cite{xu2019operational} by integrating different Markov process models for SoC valuation. We first present the storage arbitrage formulation, then describe Markov process model for two valuation models: real-time price and day-ahead price bias against real-time price (DAP-RTP bias). 
% We will demonstrate the solution method in next section.

% \subsection{Multi-stage Stochastic Arbitrage Formulation}
We formulate energy storage arbitrage using SDP as a price response problem assuming the storage can update its control over time period $t$ after observing the newest real-time market price $\lambda_t$, which the system operator will announce before the dispatch period starts. We assume the distribution if $\lambda_{t+1}$ is dependent on the price realization from the previous period $\lambda_{t}$. The arbitrage problem is formulated as
\begin{subequations}\label{eq1}
\begin{align}
    Q_{t-1}(e_{t-1}\,|\,\lambda_t) &= \max_{b_t, p_t} \lambda_t (p_t-b_t) - cp_t + V_t(e_t\,|\,\lambda_t) \\
    V_t(e_t\,|\,\lambda_t) &= \mathbb{E}_{\lambda_{t+1}}\Big[Q_{t}(e_{t}\,|\,\lambda_{t+1}) \Big| \lambda_t\Big] \label{eq:obj2}
\end{align}
where $Q_{t-1}$ is the maximized energy storage arbitrage revenue from time period $t$ till the end of the operating horizon $T$. $Q_{t-1}$ is dependent on the energy storage SoC at the end of the previous time period $e_{t-1}$, and the real-time market price $\lambda_t$. The first term in the objective function is the market income over period $t$ as the product of  $\lambda_t$ and the battery output power $(p_t-b_t)$, where $p_t$ is the discharge power and $b_t$ is the charge power. The second term of the objective function models the discharge cost, in which $c$ is the marginal discharge cost. The last term is the value-to-go function that models the opportunity value of the energy $e_t$ left in the battery. We consider a stage-dependent price model such that the future price realization $\lambda_{t+1}$ is dependent on the current price, such that $V_t$ is dependent on $\lambda_t$. \eqref{eq:obj2} shows $V_t$ is the conditional expectation of the future arbitrage revenue $Q_{t}$ given $\lambda_t$.

The arbitrage problem is subject to the following constraints
\begin{gather}
    0 \leq b_t \leq P,\; 0\leq p_t \leq P \label{p1_c2} \\
    \text{$p_t = 0$ if $\lambda_t < 0$} \label{p1_c5}\\
    \rev{e_t - e_{t-1} = -p_t/\eta^p(e_{t-1}) + b_t\eta^b(e_{t-1}) \label{p1_c1}}\\
    0 \leq e_t \leq E \label{p1_c3}
\end{gather}
\end{subequations}
where \eqref{p1_c2} models the storage power rating $P$ over charge and discharge power. \eqref{p1_c5} enforces the battery to not discharge when prices are negative, a sufficient condition for avoiding simultaneous charge and discharge~\cite{xu2019operational}. \eqref{p1_c1} models the energy storage SoC evolution constraint with \rev{discharge/charge efficiency $\eta^p$ and $\eta^b$. The efficiencies are modeled as affine functions of the starting SoC $e_{t-1}$ assuming the efficiency variation is smooth.}  \eqref{p1_c3} models the energy limit $E$ over the storage SoC.

\begin{remark}\textbf{Variable efficiency model and convexity.} Note that  \eqref{eq1} is convex with a concave-maximizing objective function despite the variable efficiency model. This is because that during each decision making stage, the starting SoC $e_{t-1}$ is not a decision variable but a given parameter,  $\eta^p(e_{t-1})$ and $\eta^b(e_{t-1})$ are parameters despite their dependency over SoC. Hence \eqref{p1_c1} is still a linear constraint and the problem is convex. SDP will provide the optimal control policy with \eqref{eq1} according to Bellman's principle of optimality~\cite{bellman1966dynamic}.
\end{remark}

% we maximize operating profit $Q_{t-1}(e_{t-1}\,|\,\lambda_t)$, $b_t$ is the energy charged into the storage over time period $t$, and $p_t$ is the discharged energy over period $t$.
% $e_t$ is the storage state-of-charge (SoC) over period $t$ modeled in \eqref{p1_c1} as the charge and discharge energy subject to the efficiency $\eta$, and \eqref{p1_c2} enforces the upper and lower energy bound. $\lambda_t$, $P$, $E$ are parameters which represent real-time price, power rating, and energy storage capacity. In this paper, we consider a generalized linear energy storage model with constant power ratings and efficiency.
% We use Markov process real-time price model and Markov process DAP-RTP bias model to calculate analytical value function $V_t(e_t|\lambda_t)$. The Markov process models and solution methods for these two models are demonstrated in following sections.

\section{Solution Method}
We introduce an analytical algorithm for solving the arbitrage problem. First, we model the stochastic price process using an order-1 Markov process and summarize how we train the Markov process using historical data. Then we reformulate the arbitrage problem based on the Markov process and present a SoC valuation method and control policy. In the end, we introduce the complete solution algorithm.

% In this section, we first illustrate training method of real-time price Markov model and DAP-RTP bias Markov model. Then a generalized solution method for arbitrage with Markov model is present, follows by variant solving schemes for different cases.

\subsection{Markov Process Model}
We model real-time price uncertainty using an order-1 Markov Process model. The Markov model includes $\mathcal{T}$ stages, and each stage has $\mathcal{N}$ nodes. We discretize real-time prices into a set of time-dependent nodes $\pi_{t,i}$, each node clusters prices within the range $[\underline{\pi}_{t,i}, \overline{\pi}_{t,i})$ where $\underline{\pi}_{t,i}$ is the lower range bound and $\overline{\pi}_{t,i}$ is the upper range bound. Note that the smallest price node does not have a lower bound, and the highest price node does not have an upper bound. Each period has the same number of nodes $\mathcal{N} = \{1,2,\dotsc,N\}$, but the node value may be different. Price nodes are associated with the following time-varying transition probability matrix $M_t = \{\rho_{i,j,t}\,|\, i\in\mathcal{N}, j\in\mathcal{N}\}$ indicating the probability that the price will transit from $\pi_{t,i}$ to $\pi_{t+1,j}$, $i \in \mathcal{N},  j \in \mathcal{N}$.

We generate the probability transition matrices $M_t$ from historical real-time price (RTP). All transition matrices in this paper are on an hourly resolution, i.e., $M_t$ are identical in same hour $H_h$, where $H_h=\{1,2,\dotsc,24\}$. This approximation is motivated by the fact that generators are scheduled hourly in day-ahead, and also to reduce the effect of limited training historical data size.
% \begin{align}
%     M_t = M_{H_h} \quad \forall t \in H_h
% \end{align}
We calculate $M_{H_h}$ base on the sum of state transit incidents happens in $H_h$ for all days in training historical data.
% Now we show how we factorize RTP data into $\mathcal{N}$ price nodes. We define negative price as $\lambda_1$ and price spike as $\lambda_N$. In this research, we set price spike as \$200/MWh, and the gap between price node is \$10/MWh, so there are 22 price nodes ($\mathcal{N} = \{1,2,\dotsc,22\}$). For example, $\lambda_t = \$15 / MWh \in \pi_{t,3}$.
% We use a subset of historical price data as our training data, usually it is historical price data prior to the testing data. 
%After grounding and factorizing RTP data into $\mathcal{N}$ state nodes, we are able to calculate hourly transition probability in each pair of states. The transition probability is given by
Given a set of training price data set $\lambda_t$, we calculate the hourly transition probability as $(\forall\, t \in H_h)$
\begin{align}\label{eq:MP}
    \rho_{i,j,H_h} = \frac{\sum_{t \in H_h}I\{\lambda_{t+1}\in\pi_{t+1,j}\}I\{ \lambda_{t}\in\pi_{t,i}\}}{\sum_{t \in H_h}I\{\lambda_{t}\in\pi_{t,i}\}}
\end{align}
where $I\{A\}$ denotes the indicator function for event $A$: $I\{A\}=1$ if $A$ is true, 0 otherwise. This transition probability is the likelihood of the next time point price in node $j$ given the current price in node $i$. Note that $\rho_{i,j,t}$ in time period $t$ is identical to transition probability of hour $H_h$ it belongs.

% \subsection{DAP-RTP Bias Markov Model Training}
\begin{remark}\textbf{DAP-RTP price bias model.}
In addition to fitting a Markov process directly over historical RTP data, we also explore Markov models using the bias between day-ahead price (DAP) and RTP, as DAP are published before the operating day starts and represent the system schedule. In this case, the training data in $\lambda_t$ in \eqref{eq:MP} becomes the historical differences between DAP and RTP, and the node $\pi_{t,i}$ are discretized based on the DAP-RTP bias. However, we solve the arbitrage problem and update the value of the price node by adding DAP over the corresponding hour to it. In this case, the node values depend on day-ahead forecast price, so the price nodes value varies for different dates. 
\end{remark}

\begin{figure}[h]
\centering
\includegraphics[width=0.8\columnwidth]{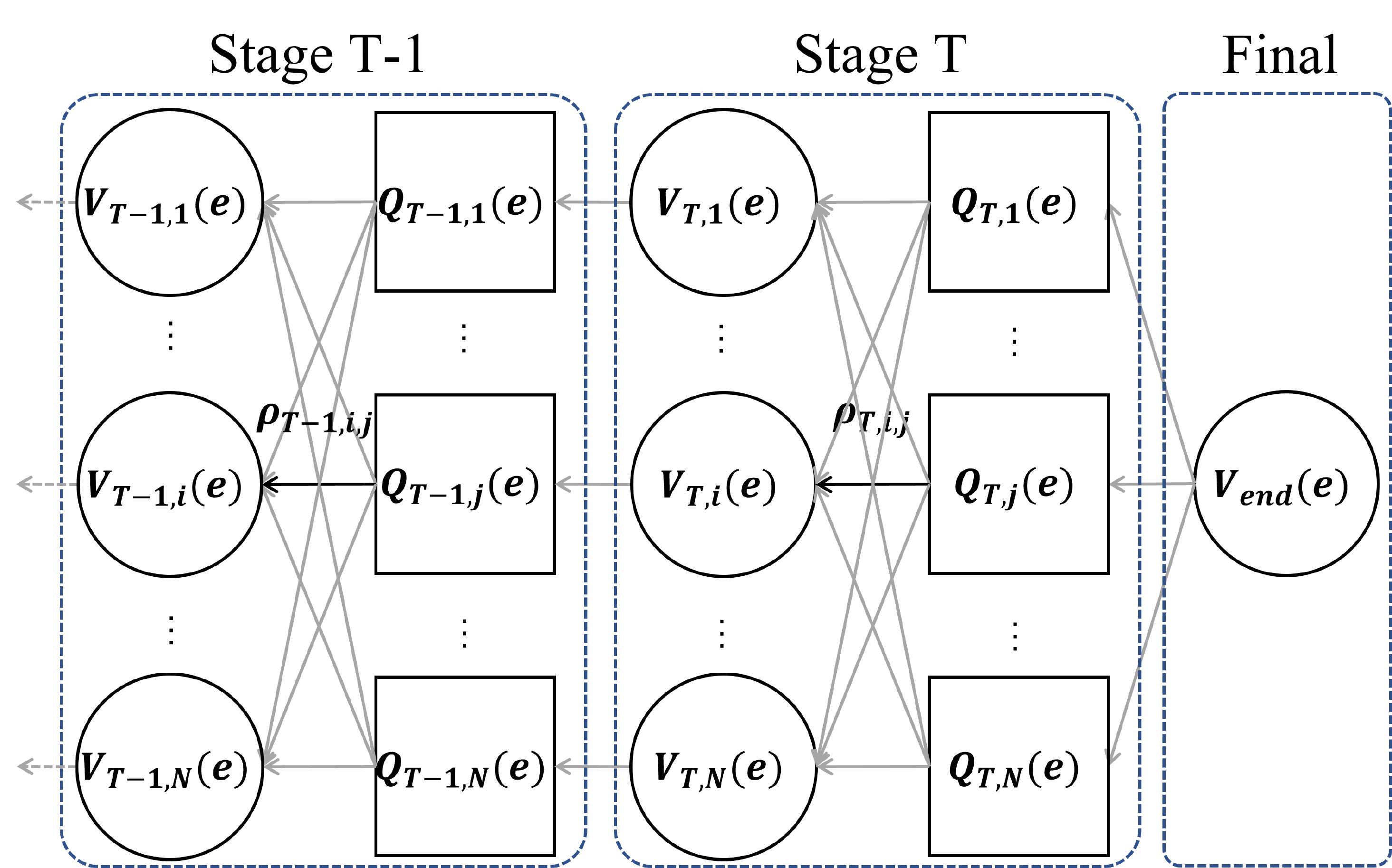}
\caption{Diagram of state transition in the proposed SDP with Markov process.}
\label{Fig.diagram}
\end{figure}

\subsection{Arbitrage Reformulation with Markov Models}
We now reformulate the price arbitrage problem in \eqref{eq1} by replacing the stochastic price $\{\lambda_t|\, t\in \mathcal{T}\}$ with the proposed Markov price process $\{\pi_{t,i}|\, t\in \mathcal{T},\, i\in \mathcal{N}\}$ as
\begin{subequations}
\begin{align}
    Q_{t-1,i}(e_{t-1}) &= \max_{b_t, p_t} \pi_{t,i}(p_t-b_t) - cp_t + V_{t,i}(e_t) \label{eq2:obj1}\\
    V_{t,i}(e_t) &= \sum_{j\in\mathcal{N}} \rho_{i,j,t}  {\cdot}Q_{t,j}(e_{t}) \label{eq2:obj2}
\end{align}\label{eq2}
\end{subequations}
subject to the same constraints \eqref{p1_c2}--\eqref{p1_c3}.
Note that the major change of the discretized formulation compared to the original one in \eqref{eq1} is that the maximized revenue $Q_{t-1,i}(e_{t-1})$ and the value-to-go function $V_{t,i}(e_t) $ now are uniquely associated with the price nodes $\pi_{t,i}$. $Q_{t-1,i}(e_{t-1})$ is the maximized arbitrage profit when the storage initial SoC is $e_{t-1}$ at price node $\pi_{t,i}$, $V_{t,i}(e_t)$ is the opportunity value of the storage SoC $e_t$ at the end of the time period $t$ at price node $\pi_{t,i}$. The conditional expectation calculation of the value function in \eqref{eq:obj2} is now replaced with \eqref{eq2:obj2} using the price node transition probability $\rho_{i,j,t}$, which models the probability of transitioning to price node $\pi_{t+1,j}$ when the current price is $\pi_{t,i}$. Hence $V_{t,i}(e_t)$ is the expected SoC value function counting in all transition possibilities associated with the current price node $\pi_{t,i}$. A simplified diagram of state transition is shown in Fig.~\ref{Fig.diagram}

\subsection{Analytical Value Function Update}

Our proposed algorithm is based on the following result that updates $q_{t-1,i}$ from $v_{t,i}$, where $q_{t-1,i}$ is the derivative of $Q_{t-1,i}$ and $v_{t,i}$ is the derivative of $V_{t,i}$
\begin{align}\label{eq3}
    &q_{t-1,i}(e) = \\
    &\begin{cases}
    v_{t,i}(e+P\eta^b)  & \text{if $\pi_{t,i}\leq v_{t,i}(e+P\eta^b){\cdot}\eta^b$} \\
    \pi_{t,i}/\eta^b  & \text{if $ v_{t,i}(e+P\eta^b){\cdot}\eta^b < \pi_{t,i} \leq v_{t,i}(e){\cdot}\eta^b$} \\
    v_{t,i}(e) & \text{if $ v_{t,i}(e){\cdot}\eta^b < \pi_{t,i} \leq [v_{t,i}(e)/\eta^p + c]^+$} \\
    (\pi_{t,i}-c)\eta^p & \text{if $ [v_{t,i}(e)/\eta^p + c]^+ < \pi_{t,i}$} \\
    & \quad\text{$ \leq [v_{t,i}(e-P/\eta^p)/\eta^p + c]^+$} \\
    v_{t,i}(e-P/\eta^p) & \text{if $\pi_{t,i} > [v_{t,i}(e-P/\eta^p)/\eta^p + c]^+$}  \nonumber
    \end{cases}
\end{align}
\rev{Note that $\eta^p$ and $\eta^b$ are dependent on the starting SoC $e$ as formulated in \eqref{eq1}, i.e., they are equivalent to $\eta^p(e)$ and $\eta^b(e)$. We omitted the function form to simplify the mathematical representation. The rest of the paper is the same on efficiency representations.} The derivation of \eqref{eq3} is summarized in the Appendix A, while an intuitive explanation is given as follows. The five listed cases in \eqref{eq3} corresponds to the five control cases in the same order: 
\begin{itemize}
    \item Storage charges at full power $P$ when price is very low;
    \item Storage charges at partial power between 0 and $P$ when price is moderately low;
    \item Storage does not charge or discharge when price is similar to the current marginal opportunity value;
    \item Storage discharges at partial power between 0 and $P$ when price is moderately high;
    \item Storage discharges at $P$ when price is very high. 
\end{itemize}
while exact conditions are as listed in \eqref{eq3}. 

To explain the principle of \eqref{eq3},  we use first-order optimality condition to solve \eqref{eq2:obj1} by taking the derivative with respect to the final SoC $e_{t-1}$ as
\begin{align}
    &q_{t-1,i}(e_{t-1}) = \nonumber\\
    &\pi_{t,i} \Big(\frac{\partial p_t}{\partial e_{t-1}}-\frac{\partial b_t}{\partial e_{t-1}}\Big) - c\frac{\partial p_t}{\partial e_{t-1}} + v_{t,i}(e_t)\frac{\partial e_t}{\partial e_{t-1}} = 0
\end{align}
and according to \eqref{p1_c1} and the Karush–Kuhn–Tucker (KKT) conditions we can conclude the following conditions
\begin{subequations}
\begin{align}
    {\partial p_t}/{\partial e_{t-1}} &= \begin{cases}
    1/\eta^p & \text{if \eqref{p1_c2} is not binding} \\
    0 & \text{if \eqref{p1_c2} is binding} \\
    \end{cases}\\
    {\partial b_t}/{\partial e_{t-1}} &= \begin{cases}
    -\eta^b & \text{if \eqref{p1_c2} is not binding} \\
    0 & \text{if \eqref{p1_c2} is binding} \\
    \end{cases}\\ 
    {\partial e_{t}}/{\partial e_{t-1}} &= \begin{cases}
    0 & \text{if \eqref{p1_c2} is not  binding} \\
    1 & \text{if \eqref{p1_c2} is binding} \\
    \end{cases}
\end{align}
\end{subequations}
now it is intuitive to see that in the second and fourth case of \eqref{eq3} the storage charge or discharge at partial power so that \eqref{p1_c2} is not binding, hence $q_{t-1,i}$ equals to either $\pi_{t,i}/\eta^b$ or $(\pi_{t,i}-c)\eta^p$ depending on the charge or discharge status; in the first and the fifth case the storage charge or discharge at full power, so the power derivatives are zero hence $q_{t-1,i}(e_{t-1}) = v_{t,i}(e_t)$, in which $e_t$ can be replaced with $e_{t-1} + P\eta^b$ or $e_{t-1}-P/\eta^p$ depending on the charge or discharge condition. In the third case both the charge and discharge power are bounded by non-negative value so $q_{t-1,i}(e_{t-1}) = v_{t,i}(e_t)$ directly holds.

\subsection{Control Policy}
After finishing the valuation calculation and obtained all value function derivatives $v_{t,i}$, we can execute the control using realized market prices $\lambda_t$. Similar to fitting the Markov price process, we first look up the corresponding value function $v_{t,i}$ associated with $\lambda_t$ by looking for the closest price node $\pi_{t,i}$ such that $\underline{\pi}_{t,i} \leq \lambda_t < \overline{\pi}_{t,i}$, then the storage control decision is updated as
\begin{subequations}

\begin{align}
    p_t &= \min\{\hat{p}_t, e_{t-1}\eta^p\} \label{eq4:1}\\
    b_t &= \min\{\hat{b}_t, (E-e_{t-1})/\eta^b\} \label{eq4:2}
\end{align}
where $\hat{p}_t$ and $\hat{b}_t$ are calculated as
\begin{align}
    &\{\hat{p}_t, \hat{b}_t\} = \nonumber\\
    &\begin{cases}
    \{0,P\} & \text{if $\lambda_t\leq v_{t,i}(e+P\eta^b)\eta^b$} \\
    \{0, \alpha\}  & \text{if $ v_{t,i}(e+P\eta)^b\eta^b < \lambda_t \leq v_{t,i}(e)\eta^b$} \\
    \{0,0\} & \text{if $ v_{t,i}(e)\eta^b < \lambda_t \leq [v_{t,i}(e)/\eta^p + c]^+$} \\
    \{\beta,0\} & \text{if $ [v_{t,i}(e)/\eta^p + c]^+ < \lambda_t$} \\
    & \quad\text{$ \leq [v_{t,i}(e-P/\eta^p)/\eta^p + c]^+$} \\
    \{P,0\} & \text{if $\lambda_t > [v_{t,i}(e-P/\eta^p)/\eta^p + c]^+$} 
    \end{cases} \label{eq4:3}
\end{align}
in which $\alpha$ and $\beta$ are given as follows
\begin{align}
    \alpha &= (v^{-1}_{t,i}(\lambda_t/\eta^b)-e_{t-1})/\eta^b \label{eq4:4}\\
    \beta &= (e_{t-1} - v^{-1}_{t,i}((\lambda_t-c)\eta^p))/\eta^p \label{eq4:5}
\end{align}
where $v^{-1}_{t,i}$ is the inverse function of $v_{t,i}$.
\end{subequations}\label{eq4}

\eqref{eq4:1} and \eqref{eq4:2} enforce the battery SoC constraints over the discharge  $\hat{p}_t$ and charge $\hat{b}_t$ decisions. \eqref{eq4:3} calculates control decisions and following the same principle as to \eqref{eq3} but use the observed price $\lambda_t$ instead of the price nodes $\pi_{t,i}$. In the first, third, and fifth cases the control decisions are clear as the storage is either at full charge power, idle, or full discharge power. In the second and fourth cases, the marginal value function $v_{t,i}$ should equal  the marginal objective value $q_{t-1,i}$, therefore we can use the inverse function of $v_{t,i}$ and the price to update the control decision.

\subsection{Solution Algorithm}
We list the complete solution algorithm while enforcing the final SoC value is higher than a given threshold $e\up{f}$. In this algorithm implementation we discretize $v_{t,i}$ into a set of $\mathcal{M}$ segments with value and SoC pairs 
\begin{align}
    \hat{v}_{t,i} = \{\nu_{t,i,m} | m\in \mathcal{M}\}
\end{align}
associated with equally spaced SoC segments $e_{t,m}$. In our implementation, we discretized the SoC into 1000 segments. The valuation algorithm is listed as following
\begin{enumerate}
    \item Set $T\to t$ to start from the last time period; initialize the final value-to-go function segments $\nu_{T,m}$ to zeros for $e_{T,m} > e\up{f}$ and to a very high value (we use \$1000/MWh) for $e_{T,m} \leq e\up{f}$. Note that the final value function does not depends on price nodes.
    \item Go to the earlier time step by setting $t-1 \to t$.
    \item During period $t$, go through each price node for $i\in\mathcal{N}$ and value function segment $m\in \mathcal{M}$. \rev{Update the charge and discharge efficiency corresponding to the SoC segment,} compute \eqref{eq3} and store $q_{t-1,i}(e)$; note that here $q_{t-1,i}(e_{t-1,m})$ is also discretized with the same granularity as the value function.
    \item Calculate the value function of the previous time step as
    \begin{align}
        \nu_{t-1,i,m} = \textstyle \sum_{j\in\mathcal{N}} \rho_{i,j,t}  q_{t-1,j}(e_{t,m})
        \label{eq13}
    \end{align}
    which is the derivative of \eqref{eq2:obj2}.
    \item Go to step 2) until reaching the first time step.
\end{enumerate}
Then after finishing the valuation, we perform the actual arbitrage simulation by updating the storage control decision using the actual system price $\lambda_t$ based on procedures described in \eqref{eq4} from the first time period until the last time period, we update the storage SoC according to \eqref{p1_c1}. 
\begin{remark}\textbf{Value function interpolation.}
It is possible to encounter price nodes with zero transition probability in specific hours on value function calculations. Therefore, we implemented a simple value interpolation with the expected value function in \eqref{eq13}. When the transition probability is all-zero, we interpolated the value function from the closest time period within the same price node.
\end{remark}

\section{Results}

We first compare different Markov process models to evaluate their performances. 
To examine the profitability robustness of the proposed method we perform arbitrage in different price zones . 
%These zones have significant difference in generation mix and transmission congestion, hence price profile and uncertainty. 
Then we trade off battery revenue with battery degradation by showing how different marginal discharge cost settings and power-to-energy ratios can impact storage arbitrage profit and the total discharged energy. \rev{We also implement arbitrage using variable efficiency energy storage, and further investigate the loss of using constant efficiency to evaluate value function. Finally, we test the model with different price and SoC granularities for the sensitivity analysis.} 

%To this end, we set 3 different P\,/\,E (1, 0.5, 0.25) and 4 different marginal discharge cost (\$0/MWh, \$10/MWh, \$30/MWh, \$50/MWh), then calculate percentage of captured profit compare to BEN-PF case with same P\,/\,E and marginal discharge cost.

\subsection{Data and Experiment Design}

\begin{figure}
\centering
\includegraphics[width=0.8\columnwidth]{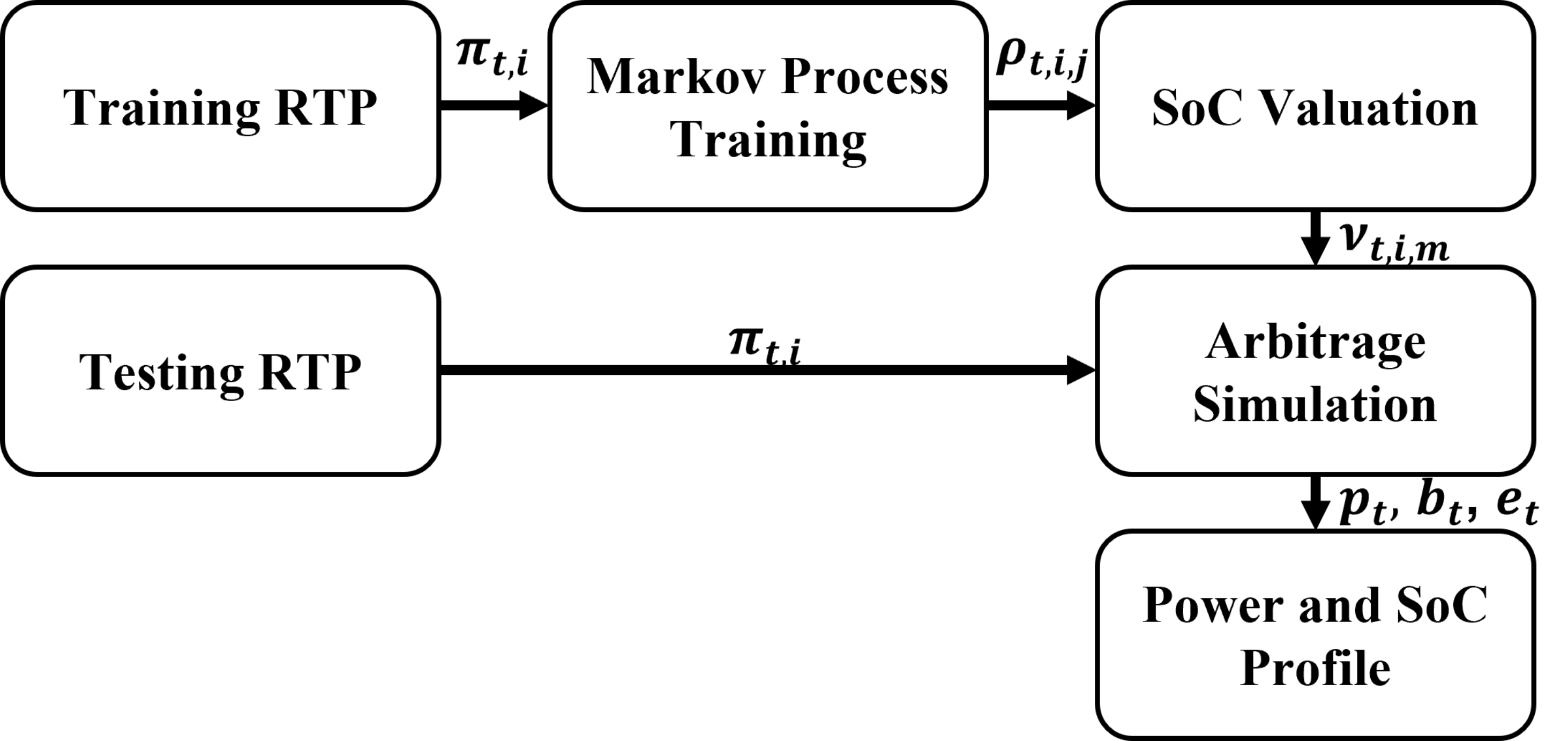}
\caption{Flow chart of the experiment set up.}
\label{Fig.flow_chart}
\end{figure}

We test the proposed arbitrage algorithm using historical price data from NYISO~\cite{nyiso}. We test the arbitrage over the 2019 data, while use data before 2019 to train the Markov process. We include prices from four zones to demonstrate arbitrage results in different price behaviors: NYC, LONGIL, NORTH, and WEST. Table~\ref{tab:zone} shows the average price and price standard deviations for the four zones in 2019. Fig.~\ref{Fig.flow_chart} shows the simulation procedure.

% The experiment follows the procedure as shown in Fig.~\ref{Fig.flow_chart}. We implement the proposed algorithm with New York ISO (NYISO) 2016-2019 historical RTP data~\cite{nyiso} in four zones: ZONE A (WEST), ZONE D (NORTH), ZONE J (NYC) and ZONE K (LONGIL). All cases in this study use 2019 RTP data as testing data. Table~\ref{tab:zone} shows the average price and price standard deviations for the four zones in 2019. Intend to investigate the effect of data size on bias, each zones are studied with different size of training data range from 1 to 3 years (i.e. [2018], [2017,2018] and [2016,2017,2018]).

\begin{table}[h]
\renewcommand{\arraystretch}{1.5}
\caption{NYISO price statistics in 2019}
\centering
% \resizebox{0.48\textwidth}{!}{
\begin{tabular}{lcccc}
\hline
\hline
Zone & NYC & LONGIL & NORTH & WEST\\
\hline
Average price                   &27.5   &32.6   &17.8   &24.7\\
Standard deviation              &28.8   &50.1   &40.2   &37.5\\
% DAP-RTP bias auto-correlation     &0.61   &0.51   &0.36   &0.50\\
% RT auto-correlation             &0.65   &0.73   &0.54   &0.55\\
% DAP-RTP correlation               &0.42   &0.30   &0.21   &0.28\\
% Num of Pos price spikes         &327	&789	&132	&578\\
% Num of Neg price spike          &208	&190	&6356	&635\\
% Num of Pos bias spike           &671	&1121	&397	&1110\\
% Num of Neg bias spike           &1061	&2570	&619	&2684\\
\hline
\hline
\end{tabular}
% }
\label{tab:zone}
\end{table}

% \begin{figure}[b]
% \centering
% \subfigure[NYC cor]{
% \includegraphics[width=0.48\columnwidth, trim = 2mm 0mm 12mm 0mm, clip]{NYC_pcor}}
% \subfigure[LONGIL cor]{
% \includegraphics[width=0.48\columnwidth, trim = 2mm 0mm 12mm 0mm, clip]{LONGIL_pcor}}
% \\
% \subfigure[NORTH cor]{
% \includegraphics[width=0.48\columnwidth, trim = 2mm 0mm 12mm 0mm, clip]{NORTH_pcor}}\label{fig:profit_north}
% \subfigure[WEST cor]{
% \includegraphics[width=0.48\columnwidth, trim = 2mm 0mm 12mm 0mm, clip]{WEST_pcor}}
% \caption{Testing and training data bias spike probability correlation in each hour}
% \label{Fig.cor}
% \end{figure}

We consider ten arbitrage settings in this study to demonstrate the impact of different stochastic models and compare them against benchmarks. The two benchmark arbitrage settings are
\begin{enumerate}
    \item \textbf{BEN-DA.} We use DAP as predictions of RTP and execute arbitrage using model predictive control. This is the first benchmark case and represents a simple arbitrage implementation.
    \item \textbf{BEN-PF.} We perform a deterministic arbitrage optimization using the actual real-time price assuming perfect forecasts. This is the second benchmark case and represents the maximum possible revenue. 
\end{enumerate}
and our proposed arbitrage algorithms
\begin{enumerate}
    \item \textbf{RT-Idp.} We perform SDP arbitrage assume there is no stage-dependency over price nodes, i.e., the transition probability does not depend on the current price node, i.e., the term $I\{\lambda_{t+1}\in\pi_{t+1,j}\}$ is removed from the Markov process training model \eqref{eq:MP}. \rev{For all RT cases, the training and testing data in Fig.~\ref{Fig.flow_chart} are real-time price.}
    \item \textbf{RT-Dep.} We perform SDP arbitrage with a Markov process trained using historical RTP data. \rev{Each stage in the Markov process has 22 price nodes, one represents price spikes, one represents negative prices, and the rest 20 are evenly spaced between \$0/MWh to \$200/MWh.}
    \item \textbf{RT-Dep-S.} Same as RT-Dep except we use different Markov processes for summer and non-summer days to represent seasonal patterns\rev{, i.e., we have two sets of $\rho_{t,i,j}$ for summer and non-summer days in Fig.~\ref{Fig.flow_chart}}
    \item \textbf{RT-Dep-W.} Same as RT-Dep except we use different Markov processes for weekdays and weekends to represent weekly patterns\rev{, i.e., we have two sets of $\rho_{t,i,j}$ for weekdays and weekends in Fig.~\ref{Fig.flow_chart}}
    \item \textbf{DB-Idp.} Similar to the RT-Idp, but the stochastic price model is trained using the DAP-RTP differences instead of using RTP directly\rev{, i.e., for all DB cases, the training and testing data are the DAP-RTP differences in Fig.~\ref{Fig.flow_chart}}  
    \item \textbf{DB-Dep.} Similar to RT-Dep,  but with a Markov process trained using historical DAP-RTP bias data. \rev{Each stage in the Markov process has 12 price nodes, with 10 nodes evenly spaced between \$-50/MWh to \$50/MWh and 2 nodes for negative and positive spikes. We use less price nodes than the RTP case because the DAP-RTP bias have a denser distribution than the RTP price.}
    \item \textbf{DB-Dep-S.} Same as DB-Dep except we use different Markov processes for summer and non-summer days to represent seasonal patterns\rev{, we have two sets of $\rho_{t,i,j}$.}
    \item \textbf{DB-Dep-W.} Same as DB-Dep except we use different Markov processes for weekdays and weekends to represent weekly patterns\rev{, we have two sets of $\rho_{t,i,j}$.}
\end{enumerate}
note that \rev{the prices used in spike nodes are the average value of the historical spikes, and mean value of the interval for other nodes,} in all DB cases we update the price node value in each day with the corresponding DAP results.

% In order to calculate the derivative of SoC value function, BEN-DA use day-ahead price forecast as real-time price in SoC valuation. In BEN-PF, we use Julia\,+\,Gurobi to solve the perfect prediction arbitrage problem. Although our algorithm also works and have similar results in BEN-PF, we want to provide an benchmark independent to our algorithm for a fair comparison. 
% For stage-independent cases (RT-Idp and DB-Idp), we use a null transition probability given by
% \begin{align}
%     \rho_{i,j,H_h} = \frac{\sum_{t \in H_h}I\{\lambda_t \in \pi_{t,j}\}}{\sum_{t \in H_h}1} \quad \forall t \in H_h
% \end{align}
% where state transit probability $\rho_{i,j,H_h}$ is independent from state node in previous time period $t$. 

\begin{remark}\textbf{Hour-ahead bidding.} Although our method targets a real-time price response problem, note that in stage-independent cases (RT-Idp and DB-Idp) we can generate all storage value functions hours before the real-time operation starts. Therefore, these two cases approximate the storage revenue potential if they participate in hour-ahead real-time market bidding. The storage bids are designed based on the value-to-go function.
\end{remark}

% Base stage-dependent cases (RT-Dep and DB-Dep) introduced in Section IV(A) .From basic dependent cases, we could obtain cases with seasonal, weekly pattern  by sub-setting original RTP/bias data to derive different state transition matrices. Each patterned cases use RTP/DAP-RTP bias data subsets to calculate SoC value function for different pattern date interval. In storage arbitrage optimization, we utilize different SoC value function according to condition of dates, i.e. whether the day is summer or weekend. The arbitrage simulation of each day (5-min resolution, 288 stage) are solved in a fraction of a second on a laptop computer for all cases investigated in this paper.

\rev{For the base battery setting, we assume 90\% one-way storage efficiency for both charge and discharge. We discuss variable efficiency energy storage in a later subsection.} 
In base cases, presumed marginal discharge cost and power-to-energy (P2E) ratio are \$10/MWh and 0.5, respectively. Note that all marginal discharge costs in this paper are not real marginal costs of specific energy storage, just presumed numbers for numerical experiments. The effect of marginal discharge cost and power-to-energy ratios will be discussed later. The Markov process training method is implemented in Python\rev{, which could be implemented within one minute per location for all different Markov process models}.
All arbitrage algorithms are implemented in Matlab except the BEN-PF case is implemented using Julia and Gurobi\rev{, the formulation of BEN-PF case could be found in Appendix B}. Although our algorithm also gives similar results in BEN-PF, we want to provide a benchmark independent of our algorithm for a fair comparison. On the computation speed, our code takes about one second on a personal computer to finish the valuation over 288 stages with 22 price nodes per stage\rev{, we further test the computation cost and model performance of one-year simulations with different price and SoC granularity in sensitivity analysis}. The Matlab implementation is available on GitHub\footnote{\url{https://github.com/niklauskun/MarkovESValuation}}.

% \begin{table}[h]
% \renewcommand{\arraystretch}{1.5}
% \caption{Cases in experiment design}
% \centering
% % \resizebox{0.48\textwidth}{!}{
% \begin{tabular}{lccc}
% \hline
% \hline
% Zone & Stage Dependency & Seasonal Pattern & Weekly Pattern\\
% \hline
% BEN-DA & F & NA & NA\\
% BEN-PF & F & NA & NA\\
% RT-Idp & F & F & F\\
% RT-Dep & T & F & F\\
% RT-Dep-S & T & T & F\\
% RT-Dep-W & T & F & T\\
% DB-Idp & F & F & F\\
% DB-Dep & T & F & F\\
% DB-Dep-S & T & T & F\\
% DB-Dep-W & T & F & T\\
% \hline
% \hline
% \end{tabular}
% % }
% \label{tab:case}
% \end{table}

\subsection{Model Comparison}
\begin{figure}[t]
\centering
\includegraphics[width=.85\columnwidth]{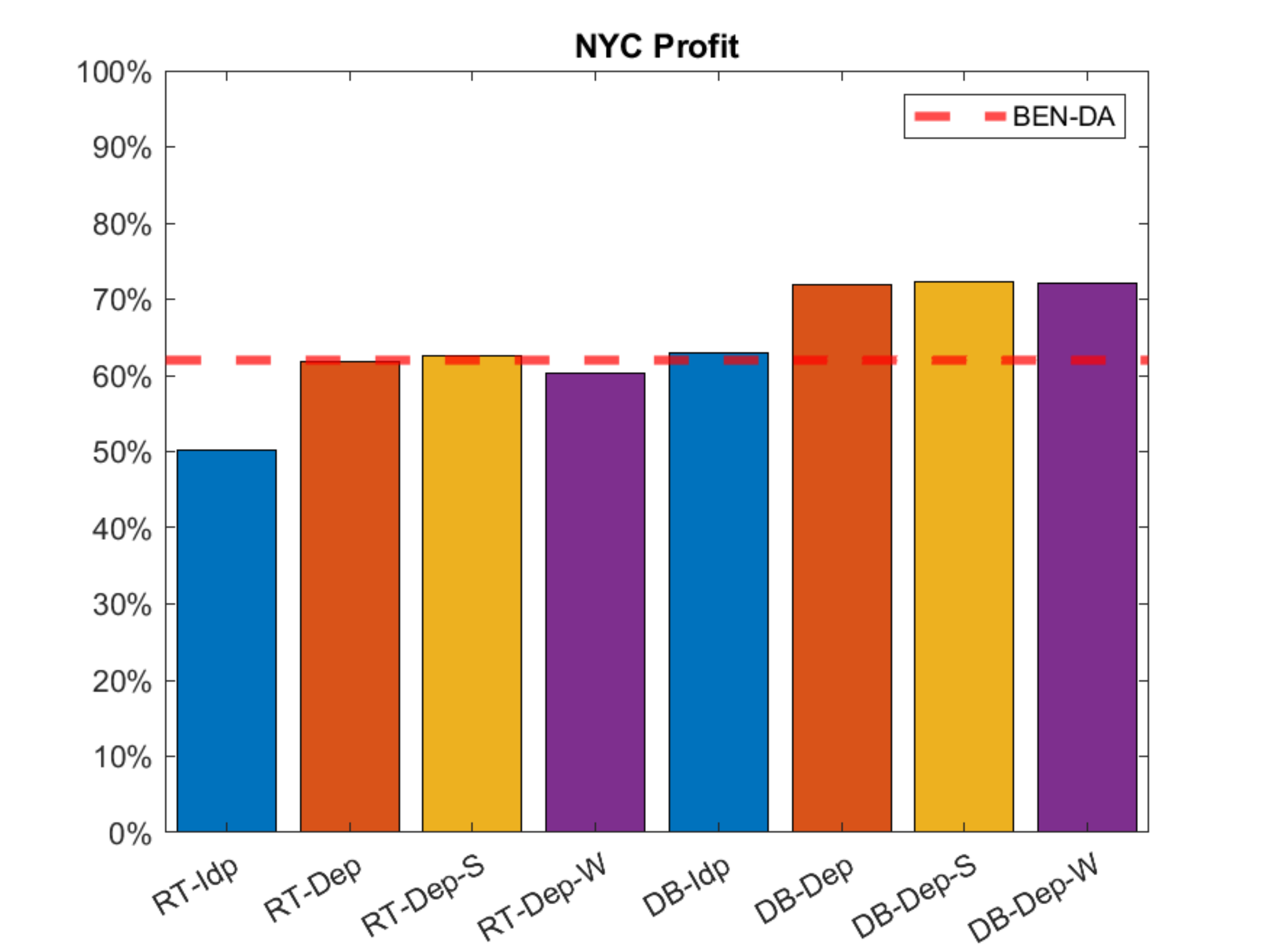}
\caption{Profit ratio captured by different cases in NYC using a 1\,MWh and 0.5 P\,/\,E battery.}
\label{Fig.ModelComparison}
\end{figure}

We show the performance of different models by comparing profits from different cases to benchmark cases (BEN-PF and BEN-DA). 
% Except specific stated, we normalize storage energy rating to 1MWh, with 0.5 P\,/\,E (2 hour storage duration), \$10/MWh presumed marginal cost of discharge, and use 2016-2019 (3 years) data as training data for all cases in this paper.
Fig.~\ref{Fig.ModelComparison} shows arbitrage profits of different models compared to BEN-PF and BEN-DA in NYC. BEN-PF is represented by 100\% in this figure as the maximum possible profit. If we consider the real-time model (RT-Dep) and the DAP-RTP bias model (DB-Dep) together, stage-dependency improves arbitrage profit significantly. Profits from all stage-dependent cases exceed BEN-DA. Profits captured by DB cases are essentially higher than BEN-DA. 

\begin{figure}[t]
\centering
\includegraphics[width=.85\columnwidth,trim = 35mm 85mm 33mm 90mm, clip]{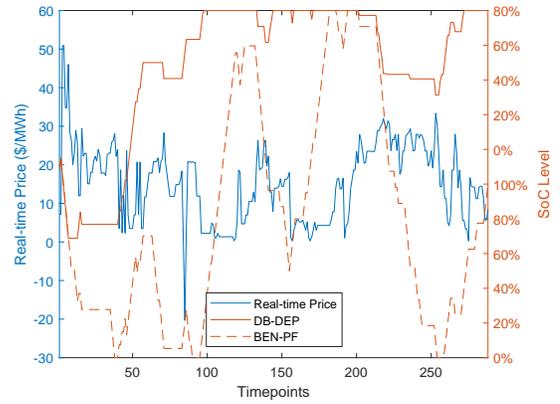}
\caption{One-day example SoC level changes of DB-Dep and BEN-PF in NYC}
\label{Fig.SOC}
\end{figure}

Seasonal pattern slightly improves the model performance with value interpolation, but it still suffers from zero transition probability \rev{in Markov process model}. The weekly pattern has a smaller subset of weekend days, considering the number of weekend days is small. Thus, it has more zero transition probabilities, and it is even detrimental to model performance. The seasonal and weekly patterns have a negligible effect on the DAP-RTP bias model, because we predict real-time prices based on day-ahead price forecast, which already captured partial seasonal and weekly patterns through different price node values in SoC valuation.

\rev{In Fig.~\ref{Fig.SOC}, we show an one-day example real-time price and SoC changes with DB-DEP and BEN-PF in NYC. Without perfect real-time price information ahead of time, DB-Dep is more conservative on discharge, therefore holds a higher SoC level compared to BEN-PF.}

In NYC, DAP-RTP bias model is superior to real-time model, but numerical result in one zone is not sufficient to generalize conclusions. In the following subsection, we test RT-Dep, DB-Idp, and DB-Dep in arbitrage simulations over different price zones, and different sizes of training data sets.

\begin{figure}[t]
\centering
\subfigure[NYC Profit]{
\includegraphics[width=0.48\columnwidth, trim = 2mm 0mm 12mm 0mm, clip]{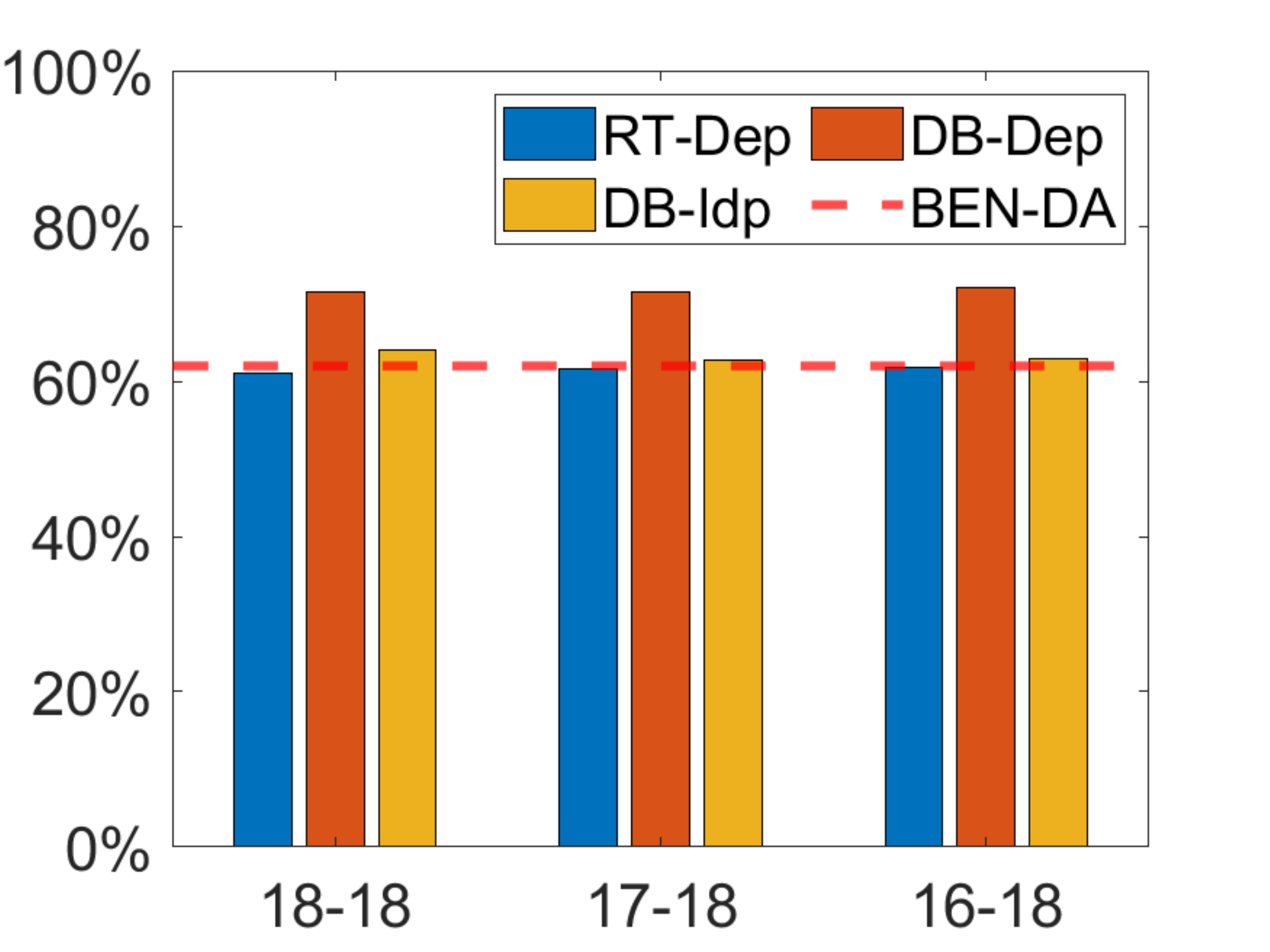}}
\subfigure[LONGIL Profit]{
\includegraphics[width=0.48\columnwidth, trim = 2mm 0mm 12mm 0mm, clip]{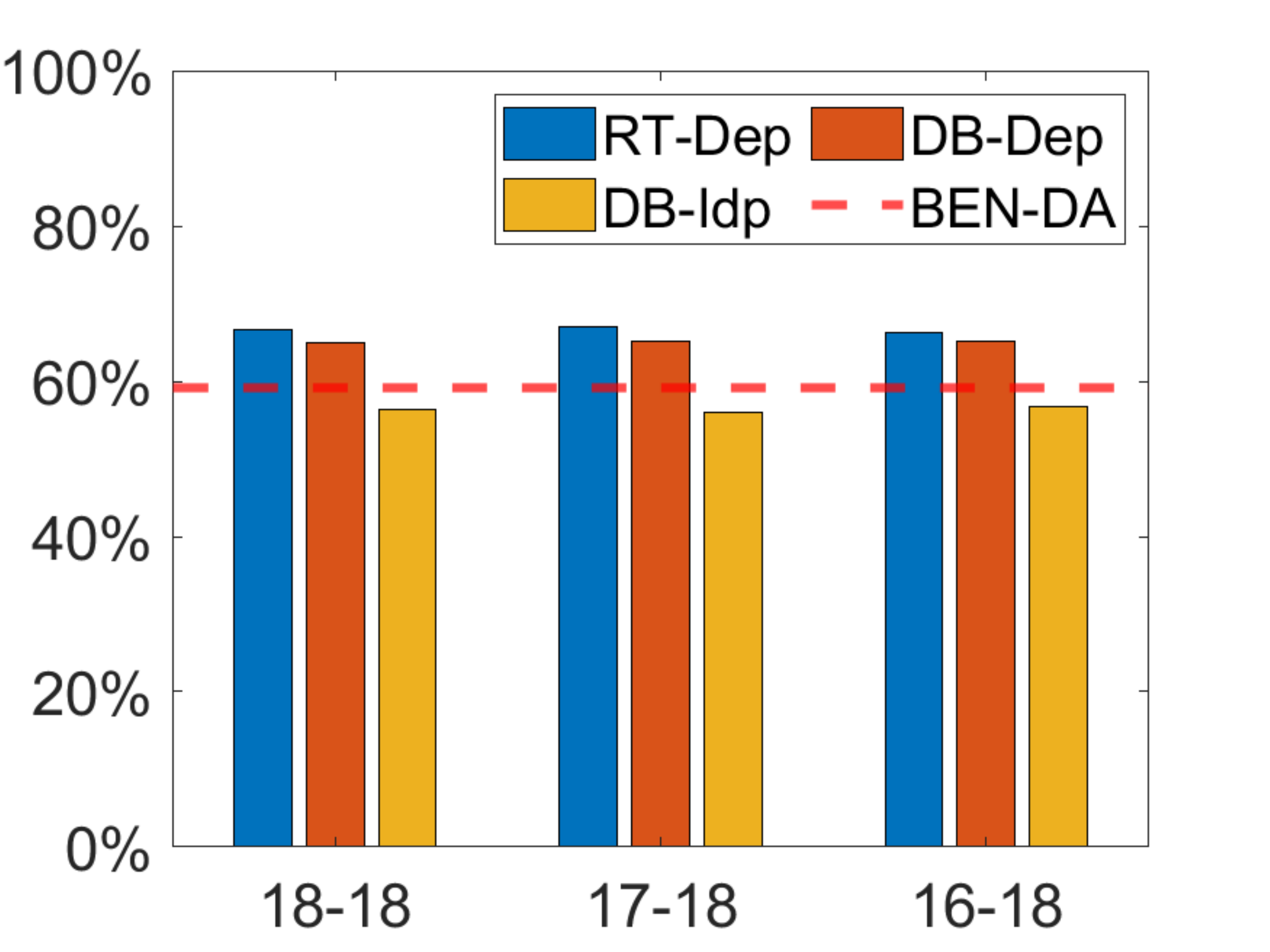}}
\\
\subfigure[NORTH Profit]{
\includegraphics[width=0.48\columnwidth, trim = 2mm 0mm 12mm 0mm, clip]{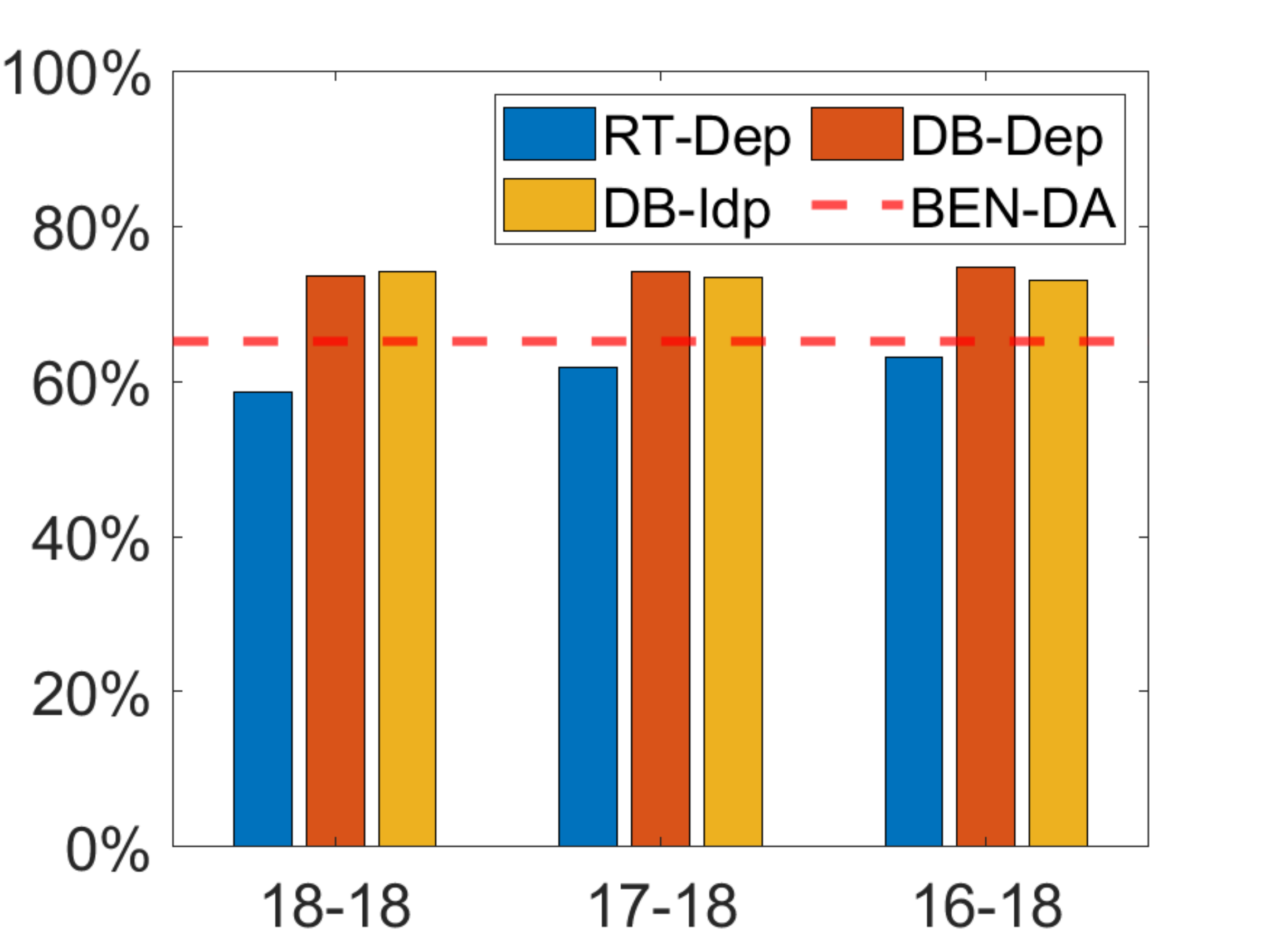}}\label{fig:profit_north}
\subfigure[WEST Profit]{
\includegraphics[width=0.48\columnwidth, trim = 2mm 0mm 12mm 0mm, clip]{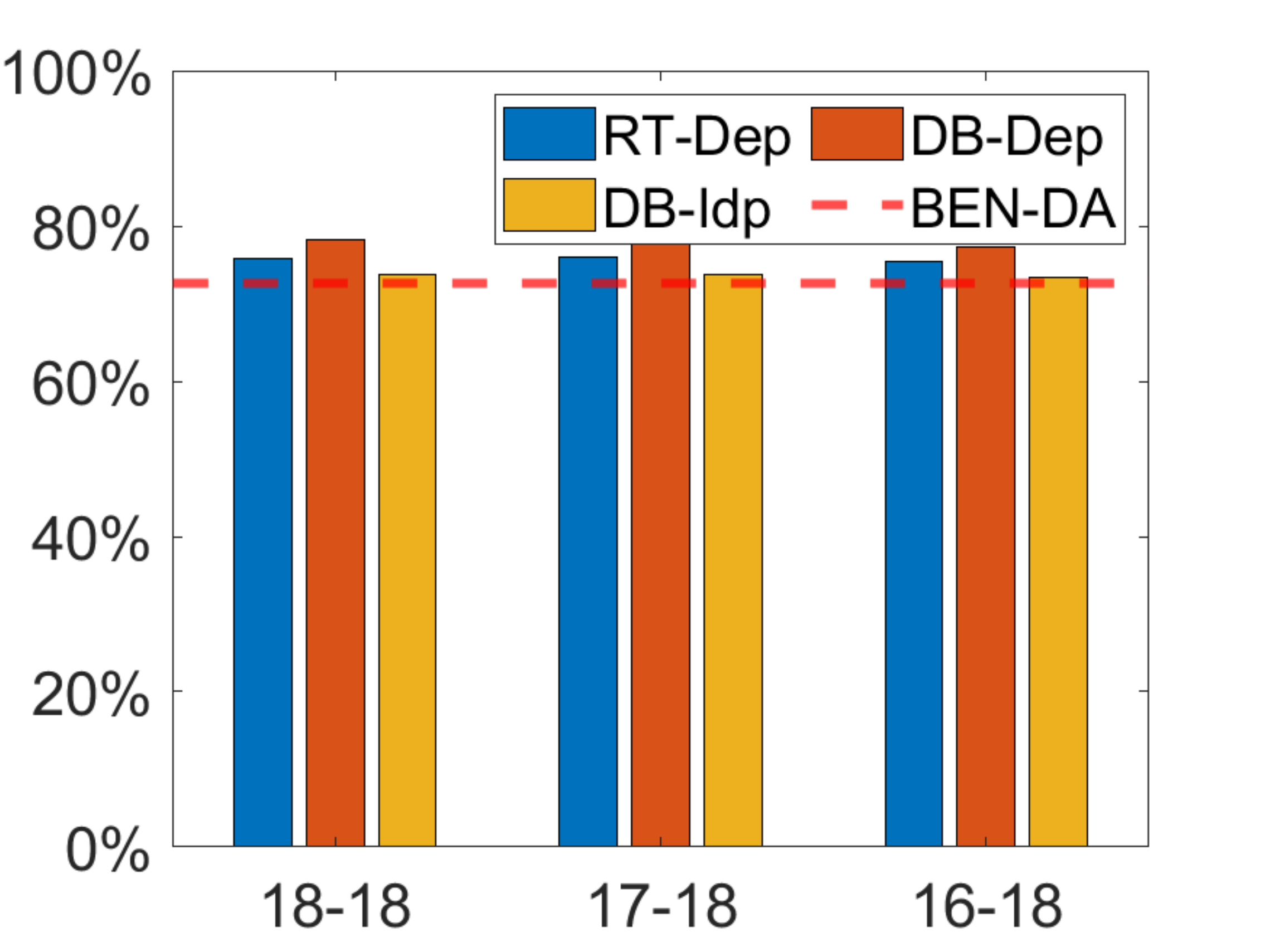}}
\caption{Profit ratio captured by real-time and day-ahead bias Markov Model, with different locations and training datasets. The training datasets are 2018 only (18-18), 2017 and 2018 (17-18), 2016 to 2018 (16-18)}
\label{Fig.Location}
\end{figure}

\begin{table*}[t]
\renewcommand{\arraystretch}{1.5}
\caption{Profit ratio, revenue and annual discharge of different power to energy ratio and presumed marginal discharge cost (\$/MWh)}
\centering
\resizebox{0.98\textwidth}{!}{
\begin{tabular}{lccccccccccccc}
\hline
\hline
\multirow{2}{*}{Zone} & \multirow{2}{*}{P2E} & \multicolumn{4}{c}{Prorated Profit Ratio [\%]} & \multicolumn{4}{c}{Revenue [k\$] / Annual Discharge [GWh]} & \multicolumn{4}{c}{Revenue per MWh discharged [\$/MWh]}\\
{}&{}&0 MC&10 MC&30 MC&50 MC&0 MC&10 MC&30 MC&50 MC&0 MC&10 MC&30 MC&50 MC\\
\hline
\multirow{3}{*}{NYC}&1&59.9&66.1&71.8&78.5&17.6\,/\,0.71 &16.8\,/\,0.25&14.1\,/\,0.11&12.9\,/\,0.07&24.91	&67.00	&125.70	&176.21\\
&0.5&67.2&72.0&78.7&84.3&11.4\,/\,0.47&10.6\,/\,0.18& 8.8\,/\,0.08&7.7\,/\,0.05&24.15	&57.71	&109.84	&158.84\\
&0.25&76.2&78.9&85.3&90.8&7.3\,/\,0.33&6.7\,/\,0.14&5.3\,/\,0.06&4.5\,/\,0.03&22.41	&47.91	&95.35	&143.71\\
\hline
\multirow{3}{*}{LONGIL}&1&56.0&59.0&62.1&62.3&28.0\,/\,0.87&27.2\,/\,0.38&23.9\,/\,0.19&21.0\,/\,0.12&32.33	&70.92	&126.83	&172.51\\
&0.5&63.5&65.1&66.7&67.4&18.8\,/\,0.58&17.8\,/\,0.28&15.0\,/\,0.13&12.8\,/\,0.08&32.21	&63.75	&117.36	&167.54
\\
&0.25&72.7&72.5&71.7&72.0&12.2\,/\,0.39&11.4\,/\,0.20&9.2\,/\,0.09&7.6\,/\,0.05&30.99	&56.75	&104.04	&154.42
\\
\hline
\multirow{3}{*}{NORTH}&1&58.4&63.5&70.1&75.3&18.3\,/\,1.12&17.1\,/\,0.34&12.6\,/\,0.09&11.0\,/\,0.04&16.36	&50.68	&142.70	&249.56
\\
&0.5&69.5&74.6&81.1&83.6&12.8\,/\,0.68&11.8\,/\,0.25&8.5\,/\,0.07&6.9\,/\,0.03&18.76	&46.74	&115.19	&200.81
\\
&0.25&79.4&83.7&90.2&88.0&8.3\,/\,0.42&7.6\,/\,0.18&5.4\,/\,0.06&4.2\,/\,0.03&19.80	&41.26	&94.18	&165.90
\\
\hline
\multirow{3}{*}{WEST}&1&67.1&70.9&75.2&78.2&36.9\,/\,1.12&35.8\,/\,0.51&30.8\,/\,0.24&27.2\,/\,0.16&32.87	&70.93	&128.12	&172.61
\\
&0.5&74.1&77.3&80.1&81.8&23.0\,/\,0.71&22.2\,/\,0.35&18.3\,/\,0.16&15.6\,/\,0.10&32.40	&63.15	&116.12	&161.22
\\
&0.25&82.2&84.2&85.8&86.7&14.2\,/\,0.45&13.5\,/\,0.24&10.9\,/\,0.11&8.9\,/\,0.06&31.72	&56.77	&102.83	&148.11
\\
\hline
\hline
\end{tabular}
}
\label{tab:profit}
\end{table*}
\subsection{Locations and Training Data Size}

\begin{figure}[h]
\centering
\includegraphics[width=0.98\columnwidth, trim = 2mm 0mm 2mm 0mm, clip]{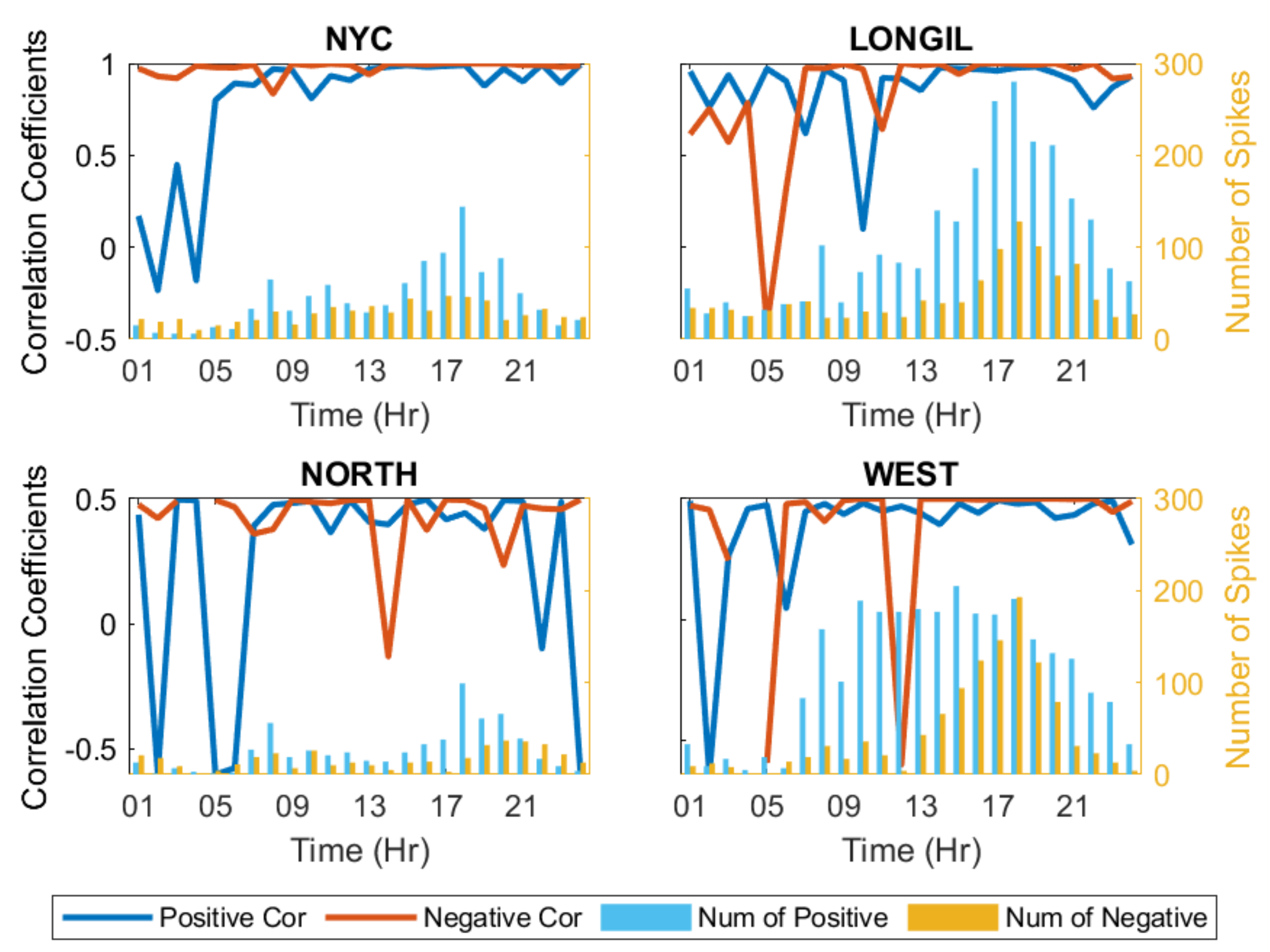}
\caption{Testing and training data bias spike number and probability correlation coefficient between testing and training data in each hour. Some correlation data not available due to 0 spikes in specific hour.}
\label{Fig.cor}
\end{figure}

We compare RT-Dep, DB-Dep, and DB-Idp cases in four different zones in NYISO, and demonstrate the impact of training data sizes. Fig.~\ref{Fig.Location} shows the percentage of captured profit compared to BEN-PF, the size of training datasets has no significant effect on model performance. The results show that one year of historical RTP data is sufficient, showing our proposed approach is also data efficient. However, suppose one wishes to capture seasonal and weekly patterns, a larger size of training data is then required to reduce incidents of zero transition probability, but we should be aware of the transformation of power systems. The training data close to testing data have more homogeneity, because the system demand, supply, and market regulation is changing over time. Yet, a larger training dataset will lead to a more stable state transition probability. Therefore, we need a trade-off between different training datasets.

Profits from stage-dependent cases are higher than BEN-DA. Notice that in Fig.~\ref{Fig.Location}(b), the model performance of RT-Dep is better than DB-Dep. The reason for this exception is LONGIL have more DAP-RTP bias spike compared to other location. Whereas in Fig.~\ref{Fig.Location}(c), NORTH has high wind penetration with lots of negative price, which make RT-Dep performance significant lower. As shown in Fig.~\ref{Fig.cor}, LONGIL and WEST have significantly more bias spike compared to other zones. However, WEST has little bias spikes in hours with low spike probability correlations, resulting in only LONGIL's DB-Dep performance being affected by the bias spike issue. We also compare DB-Dep and DB-Idp, as DB-Idp approximates hour-ahead real-time market bidding because the valuation result only depends on day-ahead forecasts. Comparing DB-Dep and DB-Idp, what stands out in all four locations is the profit loss in stage-wise dependency ignorance.

% \begin{figure}[t]
% \centering
% \subfigure[]{
% \includegraphics[width=0.9\columnwidth]{Duration MC}}
% \subfigure[]{
% \includegraphics[width=0.9\columnwidth]{Duration MC P}}
% \caption{Model performance by Different Power to Energy Ratio and Presumed Marginal Cost in NYC of a Prorated (1~MWh) Battery: (a) Comparison of DA-Dep and BEN-PF prorated profit, DA-Dep use full line, BEN-PF use dash line, (b) Ratio of profit DA-Dep captured compare to BEN-PF Profit.}
% \label{Fig.Duration}
% \end{figure}

\subsection{Power to Energy Ratio and Degradation Cost}
We investigate the effect on model performance from power-to-energy (P2E) ratios and presumed marginal discharge cost. We model three different power ratings while keeping the energy rating normalized to 1\,MWh, and include four cases with the presumed marginal discharge cost $c$ of 0, 10, 30, 50 (\$/MWh).

Table~\ref{tab:profit} compares profit ratios, revenue, and annual total discharged energy with different P2E ratios and presumed marginal discharge costs. The profit gap between our model and BEN-PF narrows down as the storage duration and marginal discharge cost increase. The results show that our proposed method can capture up to 90\% profit compared to perfect information arbitrage, with a low P2E ratio and high presumed marginal discharge cost. Table~\ref{tab:profit} also shows that arbitrage revenue per MWh increases with increment in P2E ratio and presumed marginal cost. While revenue per MWh is more sensitive to presumed marginal cost,  energy storage is prone to discharge at price spike and ignores low profitable discharge timing.

Another takeaway from the result is the comparison between revenue and the total discharged energy. Battery lifetime critically depends on energy throughput, and reducing the total energy discharged  saves  battery life significantly. For reference, currently most battery warranties require the battery to be cycled less than once per day. In our case this corresponds to keeping the annual discharged energy to less than 0.365~GWh. Our result shows that adding a \$10/MWh marginal cost can significantly reduce the energy discharged but only impact the revenue slightly. In some cases, the revenue increased after adding the \$10/MWh marginal cost as it reduced the arbitrage sensitivity to minor price fluctuations. 
By increasing the presumed marginal discharged cost from zero to \$50/MWh, we reduce the energy discharged by more than 90\%, but the revenue reduction is only  30\% to 40\%. On the other hand, the revenue per MWh discharged increased by more than ten folds. 

% Fig.~\ref{Fig.Duration} (a) compared DB-Dep profit to BEN-PF with different P\,/\,E and presumed marginal discharge cost. The profit gap between our model and BEN-PF narrow down as the storage duration and marginal discharge cost increase. Table.~\ref{tab:profit} further shows our model performance can reach around 90\% compare to perfect information arbitrage, with 0.25 P\,/\,E and \$50/MWh marginal discharge cost.  What can be clearly seen in Fig.~\ref{Fig.Duration} (b) is arbitrage revenue per MWh increase with increment in P\,/\,E ratio and presumed marginal cost. While revenue per MWh in more sensitive to presumed marginal cost, because energy storage has prone to discharge at price spike and ignore low profitable discharge timing. This can be shown in Table.~\ref{tab:profit}, the annual discharge shrinkage significantly when the marginal discharge cost in high level.

\subsection{\rev{Variable Efficiency Energy Storage}}
\begin{figure}[t]
\centering
\subfigure{
\includegraphics[width=.65\columnwidth, trim = 35mm 84mm 35mm 90mm, clip]{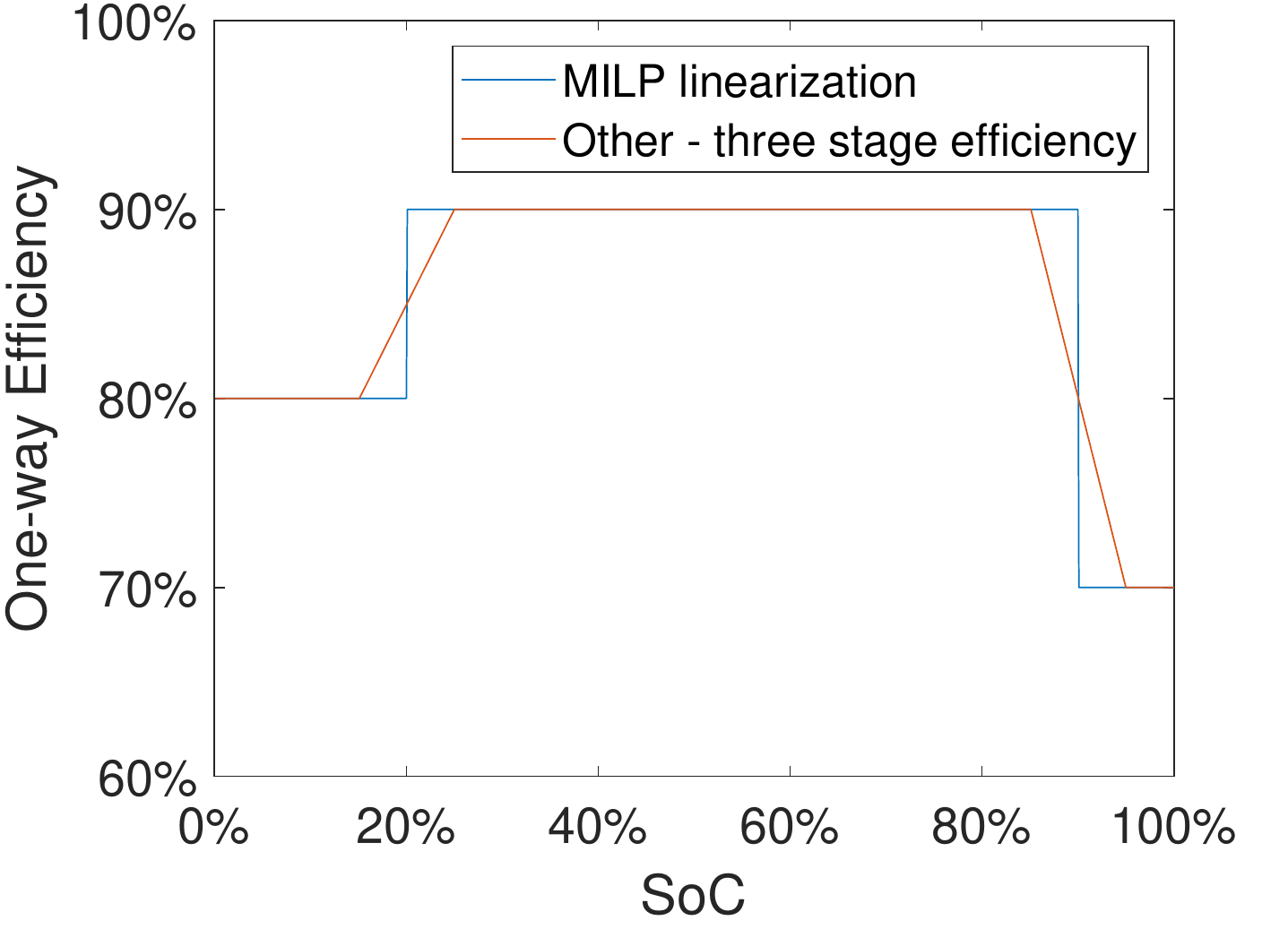}}\\
\subfigure{
\includegraphics[width=.65\columnwidth, trim = 35mm 85mm 35mm 90mm, clip]{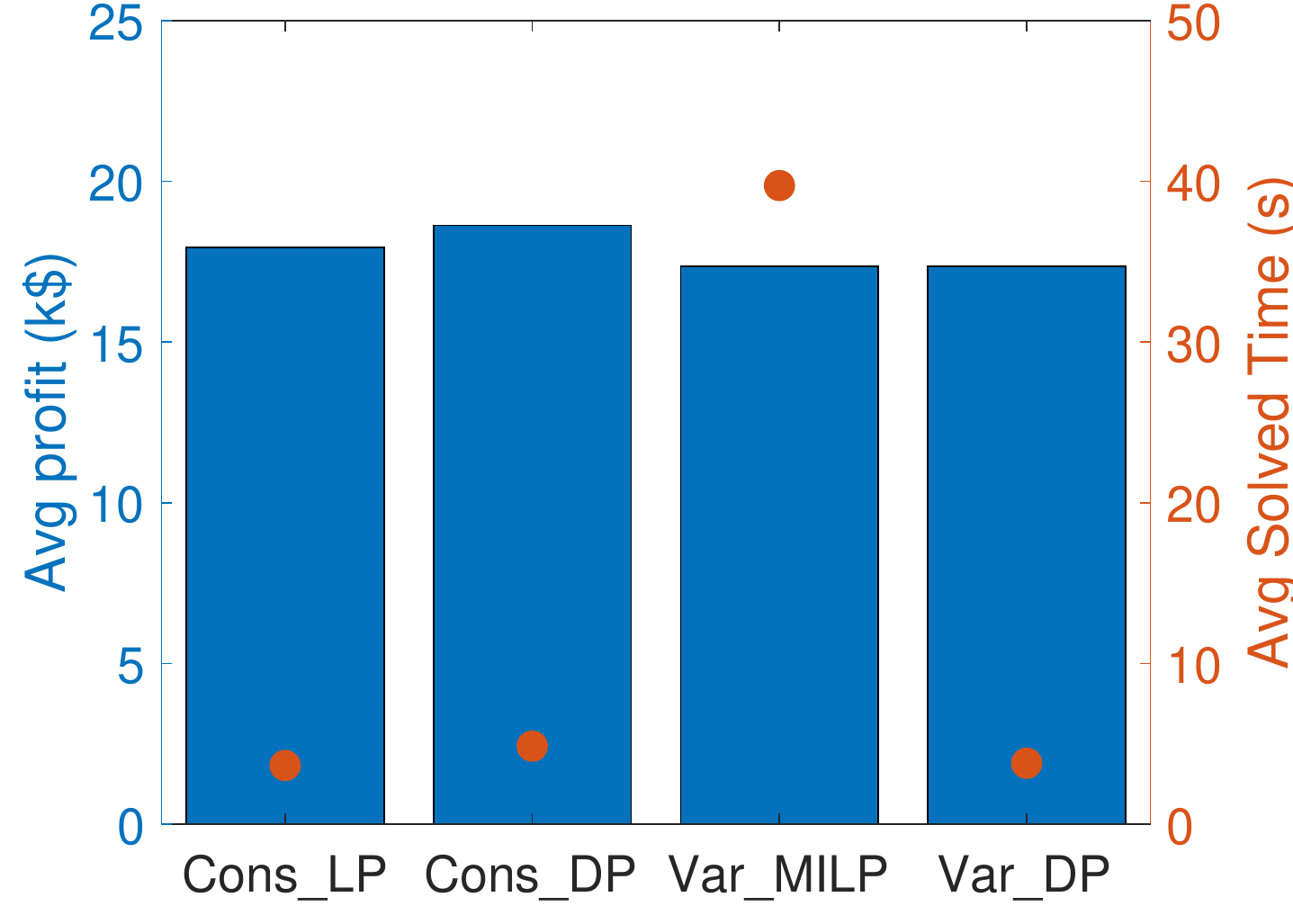}}
\caption{(a) One-way efficiency for variable efficiency case with variable efficiency and its MILP linearization  for implementation in Julia+Gurobi (b) Average deterministic arbitrage profit and solved time of all locations with LP/MILP (Julia + Gurobi) and the proposed dynamic programming (DP) algorithm. Prefix "Cons" represents constant efficiency cases; Prefix "Var" represent variable efficiency cases. Profit and solved time are shown for the entire 2019 arbitrage.}
\label{Fig. eff}
\end{figure}

\rev{We implement the DB-Dep model on a variable efficiency energy storage. In our example, we assume an energy storage has an efficiency curve dependent on SoC as shown in Fig.~\ref{Fig. eff}(a). To provide a benchmark comparison to standard mixed-integer linear program (MILP), we approximate 80\%, 90\%, 70\% one-way efficiency at 0-20\% SoC, 20-90\% SoC, and 90-100\% SoC. We present the formulation in Appendix~B.}

\rev{
We first test the two approaches using arbitrage with perfect price forecast, as MILP is limited to deterministic optimization. Fig.~\ref{Fig. eff}(b) shows arbitrage profit and computation time for the full 2019 year with perfect real-time price predictions, we get similar arbitrage profit in dynamic programming (DP) compared to MILP, which means our solution method is effective on variable efficiency energy storage. The computation times are similar in the constant efficiency case. When it comes to variable efficiency cases, solved times do not increase significantly with DP, because the algorithm and complexity are the same as constant efficiency cases. For Julia + Gurobi variable efficiency cases, we use MILP instead of LP with 3 auxiliary states, so it takes significantly longer time compared to DP or LP.} 

\begin{figure}[h]
\centering
\includegraphics[width=.65\columnwidth, trim = 35mm 85mm 35mm 90mm, clip]{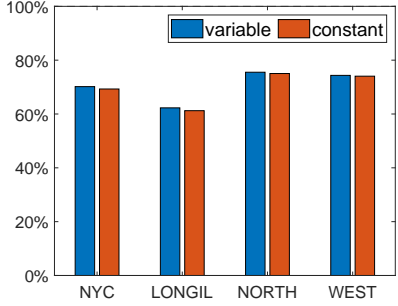}
\caption{Stochastic arbitrage profit captured compared to perfect predictions with variable efficiency storage model with variable efficiency valuation (variable) and fixed-efficiency valuation (constant).}
\label{Fig. effcomp}
\end{figure}

\rev{In order to show the distortion effect of fixing efficiency assumption on a variable efficiency energy storage, we further compared DB-Dep performance with fixed (90\% in all SoC segments) and variable efficiencies in value function evaluation, then test both methods in arbitrage with variable efficiency storage model. As shown in Fig.~\ref{Fig. effcomp}, assuming fixed efficiency slightly reduces the profit in all locations, indicating modeling variable efficiency in arbitrage will improve profitability.}

\subsection{\rev{Sensitivity Analysis}}
\begin{figure}[t]
\centering
\subfigure[Price Granularity]{
\includegraphics[width=0.65\columnwidth, trim = 35mm 83mm 35mm 90mm, clip]{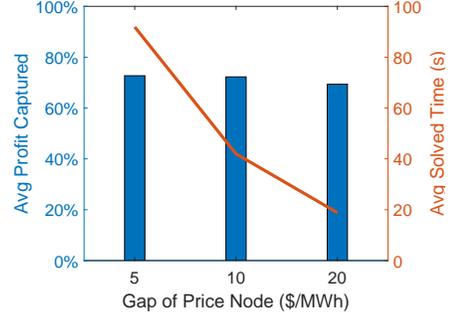}}\\
\subfigure[SoC Granularity]{
\includegraphics[width=0.65\columnwidth, trim = 35mm 83mm 35mm 90mm, clip]{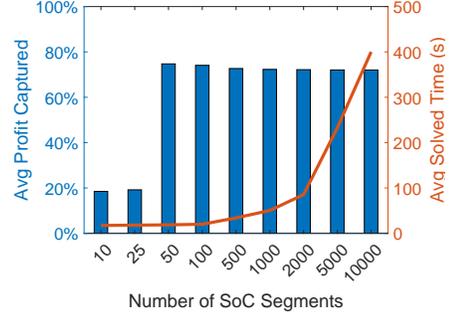}}
\caption{Sensitivity analysis of price and SoC granularity. Blue bar plot is the average profit captured compared to BEN-PF for all location, orange line plot is the average solved time for all locations. Profit and solved time are shown for the entire 2019 arbitrage. (a) price granularity sensitivity analysis, with \$5/MWh, \$10/MWh, \$20/MWh gap of price nodes; (b) SoC granularity sensitivity analysis with 9 different number of SoC segments}
\label{Fig.sa}
\end{figure}

\rev{We investigate the effect of price and SoC granularity on computational cost and model performance with DB-Dep. For price granularity, we keep upper and lower bound of bias at \$-50MWh and \$50/MWh. The gaps of each price node are \$5/MWh, \$10/MWh, and \$20/MWh, which represent 22, 12, 7 total price nodes include 2 nodes for bias spike. As shown in Fig.~\ref{Fig.sa}(a), with \$10/MWh price node gap, the computational cost decreases significantly compared to \$5/MWh gap, but maintains a close average profit ratio. Extending gap to \$20/MWh further reduces the solve times, but with a cost of profit losses. 
}

\rev{
In the SoC sensitivity analysis, we divide energy capacity into evenly spaced SoC segments. Fig.~\ref{Fig.sa}(b) shows that as the number of SoC nodes grow, computational cost increase exponentially. In this test case we used a 2-hour storage and 50 SoC nodes are sufficient to reach an optimal model performance. The optimization resolution is 5 minutes, which means in each time step the energy storage can charge/discharge 4.17\% of the total energy capacity. With 50 SoC segments, each segment represents 2\% of the energy capacity, hence the SoC granularity is fine enough for the value function evaluation. However, with longer energy storage duration, we would need a finer SoC granularity to maintain the optimal model performance. We also noted that the storage obtained slightly higher profit with 50 to 100 SoC segments compared to more segments, this is likely due to lower number of segments lead in smoother value function after piecewise linearization, and couples better with the designed Markov price process.}

\section{Conclusion}
We proposed an analytical stochastic dynamic programming algorithm considering real-time price uncertainty by integrating Markov processes. Our proposed arbitrage algorithm achieves more than 70\% profit compared to the perfect information case in the base. When considering batteries with larger energy capacity and higher degradation costs, our algorithm can reach up to 90\%  profit ratio. Our study also shows taking price auto-correlation into account can improve the profit ratio by around 10\%, while incorporating day-ahead prices into the real-time price model will provide another 10\% increase. On the other hand, seasonal and weekly patterns will not improve the performance much and require more training data.  \rev{Evaluating value function with constant efficiency can distort variable efficiency energy storage operation. According to the sensitivity analysis, both price and SoC segments granularity have an optimal threshold, which is related to the price fluctuation and energy storage setting. Over-fine granularity significantly increases the computational cost, but have little contribution to the model performance.}

% Taking account of temporal price variance such as seasonal and weekly pattern can slightly increase arbitrage profit, but creates more  zero-transition probabilities in the Markov process which  jeopardizes the model performance. Although we use interpolation method to deal with this issue, it is still not precise enough for opportunity cosy estimation in real-time arbitrage. The performance of our algorithm improve as the power to energy ratio and marginal discharge cost increase.

In the future, we plan to apply the proposed stochastic dynamic programming framework to other grid storage applications such as renewable integration, peak shaving, and stacked services. We also plan to improve the Markov process accuracy by integrating it with state-of-art price prediction models. % Finally, we would like to explore the incorporation of nonlinear battery degradation models into the proposed framework.

% Potential future directions of this work include integrating load and renewable output uncertainty, multi-stage price dependency horizon Markov process model or state-of-art price prediction model.

\section*{Appendix}

\subsection{Derivation of analytical value function update}
We start by restating result from our prior work~\cite{xu2019operational}
\begin{align}
    v_{t-1}&(e) = \mathbb{E}[q_{t-1}(e)] =\nonumber\\
    & v_t(e+P\eta^b)F_t\big(v_t(e+P\eta^b)\eta^b\big) \nonumber\\
    & + \frac{1}{\eta^b}\int_{v_t(e+P\eta^b)\eta^b}^{v_t(e)\eta^b} uf_t(u)\;du\nonumber\\
    & + v_t(e)\Big[F_t\big([v_t(e)/\eta^p + c]^+\big) -  F_t\big(v_t(e)\eta^b\big)\Big] \nonumber\\
    & + \eta^p\int_{[v_t(e)/\eta^p  +c]^+}^{[v_t(e-P/\eta^p)/\eta^p + c]^+} w f_t(w)\;dw \nonumber\\
    & - c\eta^p\Big[F_t\big({[v_t(e-P/\eta^p)}/{\eta^p} + c]^+\big)\nonumber \\
    &- F_t\big([{v_t(e)}/{\eta^p} + c]^+\big)\Big]\nonumber\\
    & + v_t(e-P/\eta^p)\Big[1-F_t\big([v_t(e-P/\eta^p)/\eta^p + c]^+\big)\Big]\label{pro1}\,
\end{align}
which calculates the value function assuming the price follows a probability distribution function $f_t$ and a cumulative distribution function $F_t$. To apply this result to the Markov process case in this paper, note that at each valuation step we are essentially solving a deterministic arbitrage at a given price node, hence $f_t$ becomes an indicator function and $F_t$ becomes a step function. \rev{For example, the first term $v_t(e+P\eta^b)F_t\big(v_t(e+P\eta^b)\eta^b\big)$ becomes $v_t(e+P\eta^b)$ if price is lower than $v_t(e+P\eta^b)\eta^b$ because $F_t\big(v_t(e+P\eta^b)\eta^b\big)$ is one only when $v_t(e+P\eta^b)\eta^b$ is larger than the price node value. Hence, repeat the derivation for all the terms in the above equation and we obtain the deterministic version used for the Markov process as shown in \eqref{eq3}.}

\subsection{MILP implementation of variable efficiency arbitrage}

\rev{
We start by showing below the multi-period arbitrage formulation, which is equivalent to our proposed stochastic dynamic programming formulation if assuming a deterministic price process $\pi_t$:}
\begin{subequations}
\begin{gather}
\max_{p_t,b_t} \quad \sum^T_t\pi_{t}{\cdot}(p_t-b_t) - cp_t
\textbf{ s.t.} \quad \text{(1c),(1e),(1f)}\nonumber
\end{gather}
\end{subequations}

\rev{We modify this model to a mixed-integer linear program model for the variable efficiency benchmark calculation with three SoC-efficiency segment pairs as}
\begin{subequations}
\begin{gather}
\max_{p_{k,t},b_{k,t}} \quad \sum^T_t\sum_k^K\pi_{t}{\cdot}(p_{k,t}-b_{k,t}) - cp_{k,t} \label{p2_obj}\\
\textbf{s.t.} \quad 0 \leq \sum_k^Kb_{k,t} \leq P,\; 0\leq \sum_k^K p_{k,t} \leq P \label{p2_c1}\\
e_{k,t} - e_{k,t-1} = -p_{k,t}/\eta^p_k + b_{k,t}\eta^b_k \label{p2_c2}\\
E_1u_{1,t} \leq e_{1,t} \leq E_1 \label{p2_c3}\\
E_ku_{k,t} \leq e_{k,t} \leq E_ku_{k-1,t}, \quad \forall k\in\{2,...,K-1\} \label{p2_c4}\\
0 \leq e_{K,t} \leq E_Ku_{K-1,t}\label{p2_c5}
\end{gather}
\end{subequations}
\rev{where $k$ is the index of different efficiency segments, we have three segments in our examples. \eqref{p2_obj} is the objective function which sums up all segments. \eqref{p2_c1} and \eqref{p2_c2} are the power rating constraints and energy storage evolution constraints implement on all three segments. \eqref{p2_c3}-\eqref{p2_c5} model the piece-wise linear efficiency model with binary variables $u_{k,t}$, which enforce the lower SoC segment must be full before upper SoC segments can take on non-zero values.}

\bibliographystyle{IEEEtran}	% (uses file "plain.bst")
\bibliography{IEEEabrv,main}		% expects file "myrefs.bib"

\begin{IEEEbiography}[{\includegraphics[width=1in,height=1.25in,clip,keepaspectratio]{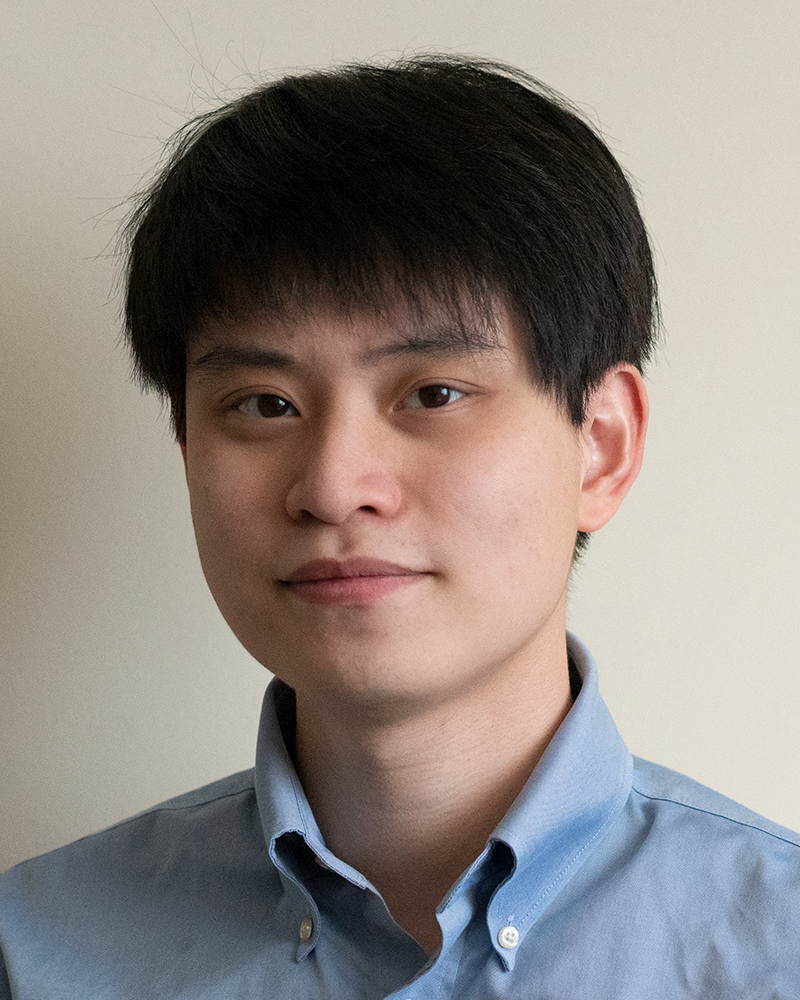}}]{Ningkun Zheng}
(S'19) received  B.S. degree in Environment and Resources Science from Zhejiang University, Zhejiang, China in 2018; M.S. degree in Environmental Health Engineering from Johns Hopkins University, Baltimore, MD, USA, in 2019. He is currently working towards the Ph.D. degree with the Department of Environmental and Earth Engineering, Columbia University, New York, NY, USA. Before joining Columbia, he was a research assistant with Carnegie Mellon Electricity Industry Center, Carnegie Mellon University, Pittsburgh, PA, USA. His research interests include power system economics, power system optimization, and energy storage.
\end{IEEEbiography}
\begin{IEEEbiography}[{\includegraphics[width=1in,height=1.25in,clip,keepaspectratio]{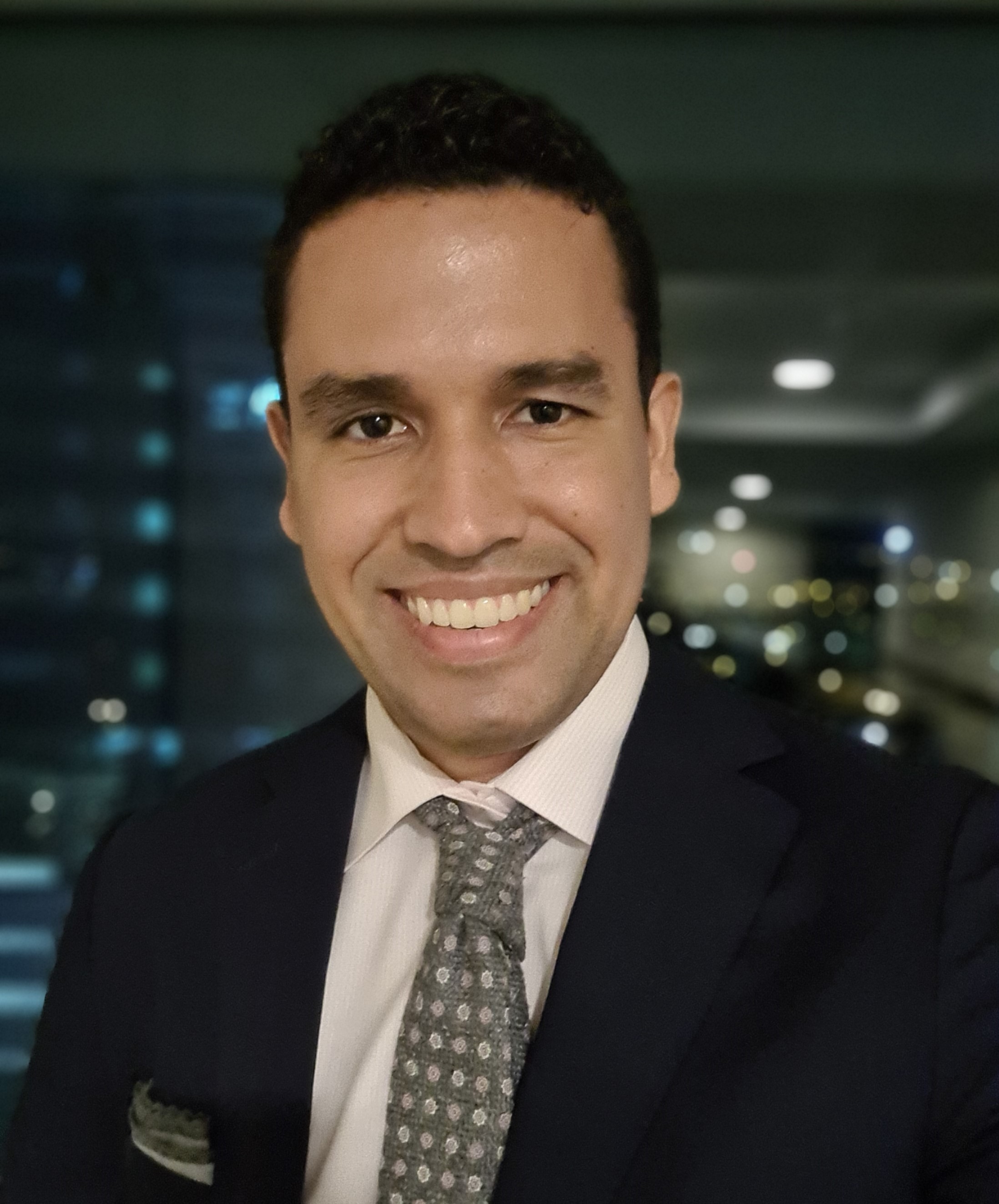}}]{Joshua Jaworski}
received his B.S. degree in Mechanical Engineering from the University of Washington, Seattle, WA, USA in 2017 and is currently working towards a M.S. degree in electrical engineering at Columbia University, New York, NY, USA. He was previously the lead process engineer at First Quantum Minerals power plant, port and copper filtration plant complex in Panama. His research interests include power systems, energy markets and distributed energy resources.
\end{IEEEbiography}
\begin{IEEEbiography}[{\includegraphics[width=1in,height=1.25in,clip,keepaspectratio]{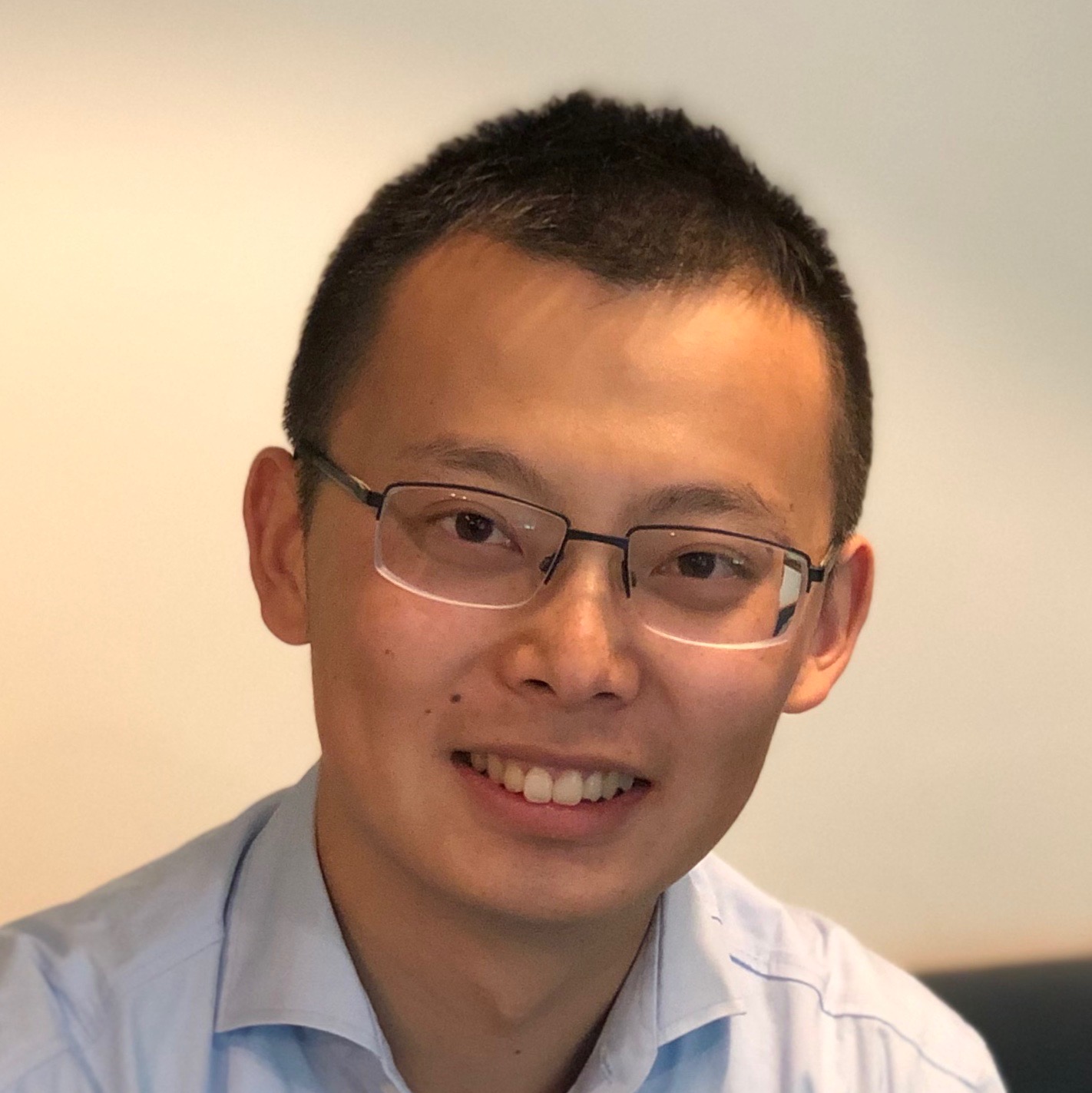}}]{Bolun Xu}
(S'14-M'18) received B.S. degrees from Shanghai Jiaotong
University, Shanghai, China in 2011;  M.Sc degree from Swiss Federal Institute of Technology, Zurich, Switzerland in 2014; and Ph.D. degree from University of Washington, Seattle, U.S. in 2018; all from Electrical Engineering.

He is currently an assistant professor in Columbia University, Department of Earth and Environmental Engineering, with affiliation in Department of Electrical Engineering. His research interests include electricity markets, energy storage, power system optimization, and power system economics. 
\end{IEEEbiography}

\end{document}